\documentclass{article}
\usepackage{amstex}
\usepackage{amssymb}
\begin{document}
\def\e#1\e{\begin{equation}#1\end{equation}}
\def\eq#1{{\rm(\ref{#1})}}
\newtheorem{thm}{Theorem}[section]
\newtheorem{lem}[thm]{Lemma}
\newtheorem{prop}[thm]{Proposition}
\newtheorem{cor}[thm]{Corollary}
\newenvironment{dfn}{\medskip\refstepcounter{thm}
\noindent{\bf Definition \thesection.\arabic{thm}\ }}{\medskip}
\newenvironment{ex}{\medskip\refstepcounter{thm}
\noindent{\bf Example \thesection.\arabic{thm}\ }}{\medskip}
\newenvironment{proof}[1][,]{\medskip\ifcat,#1
\noindent{\it Proof.\ }\else\noindent{\it Proof of #1.\ }\fi}
{\hfill$\Box$\medskip}
\def\dim{\mathop{\rm dim}}
\def\vol{\mathop{\rm vol}}
\def\Re{\mathop{\rm Re}}
\def\Im{\mathop{\rm Im}}
\def\Ker{\mathop{\rm Ker}}
\def\Aff{\mathop{\rm Aff}}
\def\Hom{\mathop{\rm Hom}}
\def\Sym{{\textstyle\mathop{\rm Sym}}}
\def\hcf{\mathop{\rm hcf}}
\def\SO{\mathop{\rm SO}}
\def\GL{\mathop{\rm GL}}
\def\U{\mathbin{\rm U}}
\def\SU{\mathop{\rm SU}}
\def\g{{\mathfrak g}}
\def\gl{\mathfrak{gl}}
\def\sn{{\textstyle\mathop{\rm sn}}}
\def\cn{{\textstyle\mathop{\rm cn}}}
\def\dn{{\textstyle\mathop{\rm dn}}}
\def\ge{\geqslant} 
\def\le{\leqslant} 
\def\cal{\mathcal}
\def\H{\mathbin{\mathbb H}}
\def\R{\mathbin{\mathbb R}}
\def\Z{\mathbin{\mathbb Z}}
\def\Q{\mathbin{\mathbb Q}}
\def\C{\mathbin{\mathbb C}}
\def\CP{\mathbb{CP}}
\def\al{\alpha}
\def\be{\beta}
\def\ga{\gamma}
\def\de{\delta}
\def\ep{\epsilon}
\def\ka{\kappa}
\def\th{\theta}
\def\la{\lambda}
\def\vp{\varphi}
\def\si{\sigma}
\def\La{\Lambda}
\def\Om{\Omega}
\def\Si{\Sigma}
\def\om{\omega}
\def\d{{\rm d}}
\def\pd{\partial}
\def\db{{\bar\partial}}
\def\ts{\textstyle}
\def\sst{\scriptscriptstyle}
\def\pha{\phantom}
\def\w{\wedge}
\def\lt{\ltimes}
\def\sm{\setminus}
\def\ov{\overline}
\def\iy{\infty}
\def\ra{\rightarrow}
\def\t{\times}
\def\ha{{\textstyle{1\over2}}}
\def\op{\oplus}
\def\ot{\otimes}
\def\bigop{\bigoplus}
\def\ti{\tilde}
\def\ms#1{\vert#1\vert^2}
\def\bms#1{\bigl\vert#1\bigr\vert^2}
\def\md#1{\vert #1 \vert}
\def\bmd#1{\bigl\vert #1 \bigr\vert}
\def\an#1{\langle #1 \rangle}
\title{Constructing special Lagrangian $m$-folds 
in $\C^m$ by evolving quadrics}
\author{Dominic Joyce, \\ Lincoln College, Oxford}
\date{}
\maketitle

\section{Introduction}
\label{qu1}

This is the second in a series of papers constructing explicit 
examples of special Lagrangian submanifolds (SL $m$-folds) in
$\C^m$. The first paper of the series \cite{Joyc2} studied
SL $m$-folds with large symmetry groups, and subsequent papers
\cite{Joyc3,Joyc4,Joyc5,Joyc6} construct examples of SL 3-folds in 
$\C^3$ using evolution equations, symmetries, ruled submanifolds
and integrable systems.

The principal motivation for these papers is to lay the foundations
for a study of the singularities of compact special Lagrangian 
$m$-folds in Calabi--Yau $m$-folds, particularly in low dimensions 
such as $m=3$. Special Lagrangian $m$-folds in $\C^m$, and especially 
{\it special Lagrangian cones}, should provide local models for 
singularities of SL $m$-folds in Calabi--Yau $m$-folds.

Understanding such singularities will be essential in making 
rigorous the explanation of Mirror Symmetry of Calabi--Yau 3-folds 
$X,\hat X$ proposed by Strominger, Yau and Zaslow \cite{SYZ}, which 
involves dual `fibrations' of $X,\hat X$ by special Lagrangian 3-tori, 
with some singular fibres. It will also be important in resolving 
conjectures made by the author \cite{Joyc1}, which attempt to define 
an invariant of Calabi--Yau 3-folds by counting special Lagrangian
homology 3-spheres.

The paper falls into three parts. The first, this section and
\S\ref{qu2}, is introductory. The second part, \S\ref{qu3}--\S\ref{qu4},
describes a general construction of special Lagrangian $m$-folds $N$ in
$\C^m$, depending on a set of {\it evolution data} $(P,\chi)$, where $P$
is an $(m\!-\!1)$-submanifold in $\R^n$. Then $N$ is the subset of $\C^m$
swept out by the image of $P$ under a 1-parameter family of linear or
affine maps $\phi_t:\R^n\ra\C^m$, which satisfy a first-order, nonlinear
o.d.e.\ in~$t$.

Examples of sets of evolution data will be given in \S\ref{qu4},
together with some progress towards a classification of such data.
The simplest interesting sets of evolution data occur when $n=m$
and $P$ is a nondegenerate quadric in $\R^m$. In the third part,
\S\ref{qu5}--\S\ref{qu7}, we apply the construction to these examples.

In this case $\phi_t(\R^m)$ must be a Lagrangian plane in $\C^m$ for
each $t$. Thus $N$ is fibred by quadrics in Lagrangian planes $\R^m$
in $\C^m$. The construction of \S\ref{qu3}--\S\ref{qu4} will also be
used in the sequel to this paper \cite{Joyc3}, with different evolution
data, to construct families of SL 3-folds in~$\C^3$.

The construction has both a linear and an affine version. In
the linear version we begin with a centred quadric $Q$ in $\R^m$,
such as an ellipsoid or a hyperboloid, and evolve its image under
linear maps $\phi_t:\R^m\ra\C^m$. This will be studied in 
\S\ref{qu5} for $\C^m$, and in more detail when $m=3$ in \S\ref{qu6}.
In the affine version we begin with a non-centred quadric $Q$ in 
$\R^m$, such as a paraboloid, and evolve its image under affine
maps $\phi_t:\R^m\ra\C^m$. This will be studied in~\S\ref{qu7}.

In some cases the family $\{\phi_t:t\in\R\}$ turns out to be 
{\it periodic} in $t$. The corresponding SL $m$-folds in $\C^m$ are 
then closed, and are interesting as local models for singular behaviour 
of SL $m$-folds in Calabi--Yau $m$-folds. Section \ref{qu55} studies 
the periodicity conditions, and proves our main result, Theorem 
\ref{qu5thm4}, on the existence of large families of SL $m$-folds 
in $\C^m$ with interesting topology, including cones on 
${\cal S}^a\t{\cal S}^b\t{\cal S}^1$ for $a+b=m-2$. When $m=3$ 
this gives many new examples of SL $T^2$-cones in $\C^3$, which 
are discussed in~\S\ref{qu6}.

In contrast to the manifolds of \cite{Joyc2}, the SL $m$-folds $N$
in $\C^m$ that we construct generically have only finite symmetry
groups. However, we shall show in \S\ref{qu43} that every set of
evolution data $(P,\chi)$ actually admits a large symmetry group
$G$, which is locally transitive on $P$. This `internal symmetry
group' does act on $N$, but not by automorphisms of $\C^m$. So we
can think of the construction as embodying a symmetry assumption,
but not of the most obvious kind.

Some of the SL $m$-folds we construct (those in \S\ref{qu5} from 
evolving ellipsoids) are already known, having been found by 
Lawlor \cite{Lawl} and completed by Harvey \cite[p.~139--143]{Harv}.
But as far as the author knows, the other examples are new. The SL 
$T^2$-cones in $\C^3$ are related to integrable systems results on 
harmonic tori in $\CP^m$. We discuss the connection in~\S\ref{qu62}.
\medskip

\noindent{\it Acknowledgements:} The author would like to thank
Nigel Hitchin, Mark Haskins, Karen Uhlenbeck, Ian McIntosh, Robert
Bryant and Chuu-Lian Terng for helpful conversations.

\section{Special Lagrangian submanifolds in $\C^m$}
\label{qu2}

We begin by defining {\it calibrations} and {\it calibrated 
submanifolds}, following Harvey and Lawson~\cite{HaLa}.

\begin{dfn} Let $(M,g)$ be a Riemannian manifold. An {\it oriented
tangent $k$-plane} $V$ on $M$ is a vector subspace $V$ of
some tangent space $T_xM$ to $M$ with $\dim V=k$, equipped
with an orientation. If $V$ is an oriented tangent $k$-plane
on $M$ then $g\vert_V$ is a Euclidean metric on $V$, so 
combining $g\vert_V$ with the orientation on $V$ gives a 
natural {\it volume form} $\vol_V$ on $V$, which is a 
$k$-form on~$V$.

Now let $\vp$ be a closed $k$-form on $M$. We say that
$\vp$ is a {\it calibration} on $M$ if for every oriented
$k$-plane $V$ on $M$ we have $\vp\vert_V\le \vol_V$. Here
$\vp\vert_V=\al\cdot\vol_V$ for some $\al\in\R$, and 
$\vp\vert_V\le\vol_V$ if $\al\le 1$. Let $N$ be an 
oriented submanifold of $M$ with dimension $k$. Then 
each tangent space $T_xN$ for $x\in N$ is an oriented
tangent $k$-plane. We say that $N$ is a {\it calibrated 
submanifold} if $\vp\vert_{T_xN}=\vol_{T_xN}$ for all~$x\in N$.
\label{qu2def1}
\end{dfn}

It is easy to show that calibrated submanifolds are automatically
{\it minimal submanifolds} \cite[Th.~II.4.2]{HaLa}. Here is the 
definition of special Lagrangian submanifolds in $\C^m$, taken
from~\cite[\S III]{HaLa}.

\begin{dfn} Let $\C^m$ have complex coordinates $(z_1,\dots,z_m)$, 
and define a metric $g$, a real 2-form $\om$ and a complex $m$-form 
$\Om$ on $\C^m$ by
\begin{align*}
g=\ms{\d z_1}+\cdots+\ms{\d z_m},\quad
\om&={i\over 2}(\d z_1\w\d\bar z_1+\cdots+\d z_m\w\d\bar z_m),\\
\text{and}\quad\Om&=\d z_1\w\cdots\w\d z_m.
\end{align*}
Then $\Re\Om$ and $\Im\Om$ are real $m$-forms on $\C^m$. Let
$L$ be an oriented real submanifold of $\C^m$ of real dimension 
$m$, and let $\th\in[0,2\pi)$. We say that $L$ is a {\it special 
Lagrangian submanifold} of $\C^m$ if $L$ is calibrated with 
respect to $\Re\Om$, in the sense of Definition \ref{qu2def1}.
We will often abbreviate `special Lagrangian' by `SL', and 
`$m$-dimensional submanifold' by `$m$-fold', so that we shall
talk about SL $m$-folds in~$\C^m$. 
\end{dfn}

As in \cite{Joyc1,Joyc2} there is also a more general definition 
of special Lagrangian submanifolds involving a {\it phase} 
${\rm e}^{i\th}$, but we will not use it in this paper. Harvey 
and Lawson \cite[Cor.~III.1.11]{HaLa} give the following 
alternative characterization of special Lagrangian submanifolds.

\begin{prop} Let\/ $L$ be a real $m$-dimensional submanifold 
of $\C^m$. Then $L$ admits an orientation making it into an
SL submanifold of\/ $\C^m$ if and only if\/ $\om\vert_L\equiv 0$ 
and\/~$\Im\Om\vert_L\equiv 0$.
\end{prop}

Note that an $m$-dimensional submanifold $L$ in $\C^m$ is 
called {\it Lagrangian} if $\om\vert_L\equiv 0$. Thus special 
Lagrangian submanifolds are Lagrangian submanifolds satisfying 
the extra condition that $\Im\Om\vert_L\equiv 0$, which is how 
they get their name.

\section{SL $m$-folds from evolution equations}
\label{qu3}

The construction of special Lagrangian $m$-folds we shall study
in this paper is based on the following theorem, which was proved 
in~\cite[Th.~3.3]{Joyc2}.

\begin{thm} Let\/ $P$ be a compact, orientable, real analytic 
$(m-1)$-manifold, $\chi$ a real analytic, nonvanishing section 
of\/ $\La^{m-1}TP$, and\/ $\phi:P\ra\C^m$ a real analytic
embedding (immersion) such that\/ $\phi^*(\om)\equiv 0$ on 
$P$. Then there exists $\ep>0$ and a unique family 
$\bigl\{\phi_t:t\in(-\ep,\ep)\bigr\}$ of real analytic maps 
$\phi_t:P\ra\C^m$ with $\phi_0=\phi$, satisfying the equation
\e
\left({\d\phi_t\over\d t}\right)^b=(\phi_t)_*(\chi)^{a_1\ldots a_{m-1}}
(\Re\Om)_{a_1\ldots a_{m-1}a_m}g^{a_mb},
\label{qu3eq1}
\e
using the index notation for (real) tensors on $\C^m$. Define 
$\Phi:(-\ep,\ep)\t P\ra\C^m$ by $\Phi(t,p)=\phi_t(p)$. Then 
$N=\Im\Phi$ is a nonsingular embedded (immersed) special 
Lagrangian submanifold of\/~$\C^m$.
\label{qu3thm1}
\end{thm}

The proof relies on a result of Harvey and Lawson 
\cite[Th.~III.5.5]{HaLa}, which says that if $P$ is a real analytic 
$(m-1)$-submanifold of $\C^m$ with $\om\vert_P\equiv 0$, then there 
is a locally unique SL submanifold $N$ containing $P$.
They assume $P$ is real analytic as their proof uses Cartan--K\"ahler
theory, which works only in the real analytic category. But this is
no loss, as by \cite[Th.~III.2.7]{HaLa} all nonsingular SL 
$m$-folds in $\C^m$ are real analytic.

We interpret equation \eq{qu3eq1} as an {\it evolution equation}
for (compact) real analytic $(m\!-\!1)$-submanifolds $\phi(P)$ of 
$\C^m$ with $\om\vert_{\phi(P)}\equiv 0$, and think of the variable $t$ 
as time. The theorem says that given such a submanifold $\phi(P)$, there 
is a 1-parameter family of diffeomorphic submanifolds $\phi_t(P)$ 
satisfying a first-order o.d.e., with $\phi_0(P)=\phi(P)$, that sweep
out an SL $m$-fold in~$\C^m$.

The condition that $P$ be compact is not always necessary in Theorem
\ref{qu3thm1}. Whether $P$ is compact or not, in a small neighbourhood 
of any $p\in P$ the maps $\phi_t$ always exist for $t\in(-\ep,\ep)$ and 
some $\ep>0$, which may depend on $p$. If $P$ is compact we can choose 
an $\ep>0$ valid for all $p$, but if $P$ is noncompact there may not 
exist such an $\ep$. If $P$ is not compact but we know for other reasons 
that there exists a family $\bigl\{\phi_t:t\in(-\ep,\ep)\bigr\}$ 
satisfying \eq{qu3eq1} and $\phi_0=\phi$, then the conclusions of 
the theorem still hold.

Now Theorem \ref{qu3thm1} should be thought of as an 
{\it infinite-dimensional} evolution problem, since the evolution 
takes place in an infinite-dimensional family of real analytic 
$(m\!-\!1)$-submanifolds. This makes the o.d.e.\ difficult to
solve explicitly, so that the theorem, in its current form, is 
unsuitable for constructing explicit SL $m$-folds. However, 
there is a method to reduce it to a {\it finite-dimensional}
evolution problem.

Suppose we can find a special class $\cal C$ of real analytic 
$(m-1)$-submanifolds $P$ of $\C^m$ with $\om\vert_P\equiv 0$, 
depending on finitely many real parameters $c_1,\ldots,c_n$, such 
that the evolution equation \eq{qu3eq1} stays within the class 
$\cal C$. Then \eq{qu3eq1} reduces to a first order o.d.e.\ on 
$c_1,\ldots,c_n$, as functions of $t$. Thus we have reduced the
infinite-dimensional problem of evolving submanifolds in $\C^m$ 
to a finite-dimensional o.d.e., which we may be able to solve 
explicitly. 

This method was used in \cite{Joyc2}, where $\cal C$ was a set of 
$(m-1)$-dimensional group orbits. We now present a more advanced 
construction based on the same idea, in which $\cal C$ consists of 
the images of an $(m-1)$-submanifold $P$ in $\R^n$ under linear or 
affine maps $\R^n\ra\C^m$. We describe the linear case first.

\begin{dfn} Let $2\le m\le n$ be integers. A set of 
{\it linear evolution data} is a pair $(P,\chi)$, where 
$P$ is an $(m-1)$-dimensional submanifold of $\R^n$, and 
$\chi:\R^n\ra\La^{m-1}\R^n$ is a linear map, such that $\chi(p)$ 
is a nonzero element of $\La^{m-1}TP$ in $\La^{m-1}\R^n$ for each 
nonsingular point $p\in P$. We suppose also that $P$ is not 
contained in any proper vector subspace $\R^k$ of~$\R^n$.

Let $\Hom(\R^n,\C^m)$ be the real vector space of linear maps
$\phi:\R^n\ra\C^m$, and define ${\cal C}_P$ to be the subset of 
$\phi\in\Hom(\R^n,\C^m)$ such that 
\begin{itemize}
\item[(i)] $\phi^*(\om)\vert_P\equiv 0$, and
\item[(ii)] $\phi\vert_{T_pP}:T_pP\ra\C^m$ is injective 
for all $p$ in a dense open subset of~$P$.
\end{itemize}

If $\phi\in\Hom(\R^n,\C^m)$ then (i) holds if and only if 
$\phi^*(\om)\in V_P$, where $V_P$ is the vector subspace of 
elements of $\La^2(\R^n)^*$ which restrict to zero on $P$. This 
is a quadratic condition on $\phi$. Also (ii) is an open condition
on $\phi$. Thus ${\cal C}_P$ is an open set in the intersection of 
a finite number of quadrics in $\Hom(\R^n,\C^m)$. Let $\R^m$ be a 
Lagrangian plane in $\C^m$. Then any linear map $\phi:\R^n\ra\R^m$ 
satisfies (i), and generic linear maps $\phi:\R^n\ra\R^m$ satisfy 
(ii). Hence ${\cal C}_P$ is nonempty.
\label{qu3def1}
\end{dfn}

Note that the requirement that $\chi$ be both linear in $\R^n$
and tangent to $P$ at every point is a very strong condition
on $P$ and $\chi$. Thus sets of linear evolution data are quite
rigid things, and not that easy to construct. We will give some
examples in \S\ref{qu4}. First we show how to construct SL
$m$-folds in $\C^m$ using linear evolution data.

\begin{thm} Let\/ $(P,\chi)$ be a set of linear evolution data, 
and use the notation above. Suppose $\phi\in{\cal C}_P$. Then 
there exists $\ep>0$ and a unique real analytic family 
$\bigl\{\phi_t:t\in(-\ep,\ep)\bigr\}$ in ${\cal C}_P$ with\/ 
$\phi_0=\phi$, satisfying the equation
\e
\left({\d\phi_t\over\d t}(x)\right)^b=(\phi_t)_*(\chi(x))^{a_1\ldots a_{m-1}}
(\Re\Om)_{a_1\ldots a_{m-1}a_m}g^{a_mb}
\label{qu3eq2}
\e
for all\/ $x\in\R^n$, using the index notation for tensors in $\C^m$.
Furthermore, $N=\bigl\{\phi_t(p):t\in(-\ep,\ep)$, $p\in P\bigr\}$ is 
a special Lagrangian submanifold in $\C^m$ wherever it is nonsingular. 
\label{qu3thm2}
\end{thm}

Before we prove the theorem, here are some remarks about it.
Equation \eq{qu3eq2} is a first-order o.d.e.\ upon $\phi_t$, and 
should be compared with equation \eq{qu3eq1} of Theorem \ref{qu3thm1}.
The key point to note is that as $\chi$ is linear, the right 
hand side of \eq{qu3eq2} is linear in $x$, and so \eq{qu3eq2}
makes sense as an evolution equation for linear maps $\phi_t$.
However, the right hand side of \eq{qu3eq2} is a homogeneous
polynomial of order $m-1$ in $\phi_t$, so for $m>2$ it is
a {\it nonlinear}~o.d.e.

Also observe that \eq{qu3eq2} works for $\phi_t$ in $\Hom(\R^n,\C^m)$, 
and not just ${\cal C}_P$. If the evolution starts in ${\cal C}_P$, 
then it stays in ${\cal C}_P$ for small $t$. But it can be helpful to 
think of the evolution as happening in $\Hom(\R^n,\C^m)$ rather than in 
${\cal C}_P$, because ${\cal C}_P$ may be singular, but $\Hom(\R^n,\C^m)$ 
is nonsingular. Thus, we do not run into problems when the evolution hits 
a singular point of~${\cal C}_P$.

\begin{proof}[Theorem \ref{qu3thm2}] As above, equation \eq{qu3eq2} is 
a well-defined, first-order o.d.e.\ upon $\phi_t$ in $\Hom(\R^n,\C^m)$ 
of the form ${\d\phi_t\over\d t}=Q(\phi_t)$, where $Q:\Hom(\R^n,\C^m)
\ra\Hom(\R^n,\C^m)$ is a homogeneous polynomial of degree $m-1$. The
existence for some $\ep>0$ of a unique, real analytic solution 
$\bigl\{\phi_t:t\in(-\ep,\ep)\bigr\}$ in $\Hom(\R^n,\C^m)$ with 
initial value $\phi_0=\phi$ follows easily from standard results 
on ordinary differential equations.

The rest of the proof follows that of Theorem \ref{qu3thm1}, given
in \cite[Th.~3.3]{Joyc2}, with small modifications. The compactness 
of $P$ in Theorem \ref{qu3thm1} was used only to prove existence of 
the family $\bigl\{\phi_t:t\in(-\ep,\ep)\bigr\}$, which we have 
already established, so we don't need to suppose $P$ is compact. 
The evolution equation \eq{qu3eq1} in Theorem \ref{qu3thm1} is 
exactly the restriction of \eq{qu3eq2} from $\R^n$ to $P$. Thus the 
proof in Theorem \ref{qu3thm1} that $N$ is special Lagrangian also 
applies here, wherever $N$ is nonsingular.

It remains only to show that $\bigl\{\phi_t:t\in(-\ep,\ep)\bigr\}$ 
lies in ${\cal C}_P$, rather than just in $\Hom(\R^n,\C^m)$. Now
$\om\vert_N\equiv 0$ as $N$ is special Lagrangian, and this
implies that $\phi_t^*(\om)\vert_P\equiv 0$ for $t\in(-\ep,\ep)$.
So part (i) of Definition \ref{qu3def1} holds for $\phi_t$. But
part (ii) is an open condition, and it holds for $\phi_0=\phi$ as 
$\phi\in{\cal C}_P$. Thus, making $\ep>0$ smaller if necessary,
we see that $\phi_t\in{\cal C}_P$ for all~$t\in(-\ep,\ep)$.
\end{proof}

Next we generalize the ideas above from linear to {\it affine} 
(linear+constant) maps $\phi$. Here are the analogues of Definition 
\ref{qu3def1} and Theorem~\ref{qu3thm2}.

\begin{dfn} Let $2\le m\le n$ be integers. A set of 
{\it affine evolution data} is a pair $(P,\chi)$, where 
$P$ is an $(m-1)$-dimensional submanifold of $\R^n$, and 
$\chi:\R^n\ra\La^{m-1}\R^n$ is an affine map, such that 
$\chi(p)$ is a nonzero element of $\La^{m-1}TP$ in 
$\La^{m-1}\R^n$ for each nonsingular $p\in P$. We suppose 
also that $P$ is not contained in any proper affine subspace 
$\R^k$ of~$\R^n$.

Let $\Aff(\R^n,\C^m)$ be the affine space of affine maps 
$\phi:\R^n\ra\C^m$, and define ${\cal C}_P$ to be the subset 
of $\phi\in\Aff(\R^n,\C^m)$ satisfying parts (i) and (ii) of
Definition \ref{qu3def1}. Then ${\cal C}_P$ is nonempty, 
and is an open set in the intersection of a finite number 
of quadrics in~$\Aff(\R^n,\C^m)$. 
\label{qu3def2}
\end{dfn}

\begin{thm} Let\/ $(P,\chi)$ be a set of affine evolution data, 
and use the notation above. Suppose $\phi\in{\cal C}_P$. Then 
there exists $\ep>0$ and a unique real analytic family 
$\bigl\{\phi_t:t\in(-\ep,\ep)\bigr\}$ in ${\cal C}_P$ with\/ 
$\phi_0=\phi$, satisfying \eq{qu3eq2} for all\/ $x\in\R^n$, 
using the index notation for tensors in $\C^m$. Furthermore, 
$N=\bigl\{\phi_t(p):t\in(-\ep,\ep)$, $p\in P\bigr\}$ is a 
special Lagrangian submanifold in $\C^m$ wherever it is nonsingular. 
\label{qu3thm3}
\end{thm}

Now the affine case in $\R^n$ can in fact be reduced to the linear 
case in $\R^{n+1}$, by regarding $\R^n$ as the hyperplane $\R^n\t\{1\}$
in $\R^{n+1}=\R^n\t\R$. Then any affine map $\phi:\R^n\ra\C^m$ extends to 
a unique linear map $\phi':\R^{n+1}\ra\C^m$. Thus Theorem \ref{qu3thm3}
follows immediately from Theorem~\ref{qu3thm2}.

\section{Examples of evolution data}
\label{qu4}

We now give examples of sets of linear and affine evolution data
$(P,\chi)$, in order to apply the construction of \S\ref{qu3}.
We begin in \S\ref{qu41} by showing that quadrics in $\R^m$ are 
examples of evolution data with $m=n$. The corresponding SL 
$m$-folds will be studied in~\S\ref{qu5}--\S\ref{qu7}. 

Section \ref{qu42} gives two trivial examples of evolution data, 
and classifies sets of evolution data in the cases $m=2$ and $m=n$. 
Then \S\ref{qu43} considers the {\it symmetries} of sets of evolution 
data, and shows that every set of evolution data $(P,\chi)$ has a 
large symmetry group $G$ which acts locally transitively on $P$. 
Finally \S\ref{qu44} discusses the classification of evolution
data, and the r\^ole of the symmetry group.

\subsection{Quadrics in $\R^m$ as examples of evolution data}
\label{qu41}

A large class of examples of evolution data arise as {\it quadrics} 
in $\R^m$, with~$n=m$. 

\begin{thm} Let\/ $\R^m$ have coordinates $(x_1,\ldots,x_m)$,
and for $j=1,\ldots,m$ define $e_j\in\R^m$ by $x_j=1$ and\/ 
$x_k=0$ for $j\ne k$. Let\/ $Q:\R^m\ra\R$ be a quadratic 
polynomial. Define $\chi:\R^m\ra\La^{m-1}\R^m$ by
\e
\begin{split}
\!\!\!\!\!\chi(x)&=\d Q(x)\cdot(e_1\w\cdots\w e_m)\\
&=\sum_{j=1}^m(-1)^{j-1}{\pd Q(x)\over\pd x_j}e_1\w\cdots\w
e_{j-1}\w e_{j+1}\w\cdots\w e_m.
\end{split}
\label{qu4eq1}
\e
Let\/ $P$ be the quadric $\bigl\{x\in\R^n:Q(x)=c\bigr\}$ for some
$c\in\R$, and suppose $P$ is nonempty and nondegenerate.

If\/ $Q$ is a homogeneous quadratic polynomial then $(P,\chi)$ 
is a set of linear evolution data in the sense of Definition
\ref{qu3def1} with\/ $n=m$, and otherwise $(P,\chi)$ is a set of 
affine evolution data in the sense of Definition \ref{qu3def2}
with\/~$n=m$.
\label{qu4thm1}
\end{thm}

The proof of this theorem is simple. As $Q$ is quadratic, 
$\d Q$ is linear or affine, so $\chi(x)$ is linear or affine 
in $x$. Since $\chi=\d Q\cdot(e_1\w\cdots\w e_m)$ and $P$ is a 
level set of $Q$, it is clear that $\chi$ lies in $\La^{m-1}TP$ 
on $P$. We leave the details to the reader. 

Here are three examples in $\R^m$, using notation as above.

\begin{ex} Let $1\le a\le m$, and define $P$ and $\chi$ by
\begin{align*}
P&=\bigl\{(x_1,\ldots,x_m)\in\R^m:
x_1^2+\cdots+x_a^2-x_{a+1}^2-\cdots-x_m^2=1\bigr\},\\
\chi&=2\sum_{j=1}^a(-1)^{j-1}x_j\,e_1\w\cdots\w e_{j-1}\w 
e_{j+1}\w\cdots\w e_m\\
&-2\sum_{j=a+1}^m(-1)^{j-1}x_j\,e_1\w\cdots\w e_{j-1}\w 
e_{j+1}\w\cdots\w e_m.
\end{align*}
Then $P$ is nonsingular in $\R^m$, and $(P,\chi)$ is a set of 
linear evolution data.
\label{qu4ex1}
\end{ex}

\begin{ex} Let $m/2\le a<m$, and define $P$ and $\chi$ by
\begin{align*}
P&=\bigl\{(x_1,\ldots,x_m)\in\R^m:
x_1^2+\cdots+x_a^2-x_{a+1}^2-\cdots-x_m^2=0\bigr\},\\
\chi&=2\sum_{j=1}^a(-1)^{j-1}x_j\,e_1\w\cdots\w e_{j-1}\w 
e_{j+1}\w\cdots\w e_m\\
&-2\sum_{j=a+1}^m(-1)^{j-1}x_j\,e_1\w\cdots\w e_{j-1}\w 
e_{j+1}\w\cdots\w e_m.
\end{align*}
Then $P$ is a quadric cone in $\R^m$ with an isolated singular point 
at 0, and $(P,\chi)$ is a set of linear evolution data.
\label{qu4ex2}
\end{ex}

\begin{ex} Let $(m-1)/2\le a\le m-1$, and define $P$ and $\chi$ by
\begin{align*}
P&=\bigl\{(x_1,\ldots,x_m)\in\R^m:
x_1^2+\cdots+x_a^2-x_{a+1}^2-\cdots-x_{m-1}^2+2x_m=0\bigr\},\\
\chi&=2(-1)^{m-1}e_1\!\w\!\cdots\!\w\!e_{m-1}\!+\!
2\sum_{j=1}^a(-1)^{j-1}x_je_1\!\w\!\cdots\!\w\!e_{j-1}\!\w\!
e_{j+1}\!\w\!\cdots\!\w\!e_m\\
&\qquad\qquad-2\sum_{j=a+1}^{m-1}(-1)^{j-1}x_j\,e_1\w\cdots\w 
e_{j-1}\w e_{j+1}\w\cdots\w e_m.
\end{align*}
Then $P$ is nonsingular in $\R^m$, and $(P,\chi)$ is a set of 
affine evolution data.
\label{qu4ex3}
\end{ex}

The classifications of centred quadrics in $\R^m$ up to linear
automorphisms, and of general quadrics in $\R^m$ up to affine
automorphisms, are well known. Our construction is unchanged  
under linear or affine automorphisms of $\R^m$. It can be 
shown that all interesting sets of evolution data arising
from Theorem \ref{qu4thm1} are isomorphic to one of the 
cases of Examples \ref{qu4ex1}--\ref{qu4ex3}, under an 
affine automorphism of $\R^m$ and a rescaling of~$\chi$.

Here we exclude quadrics admitting a translational symmetry
group $\R^k$ for $k\ge 1$ as uninteresting, since they lead
to special Lagrangian submanifolds $N$ in $\C^m$ with the same 
translational symmetry group. It then follows that $N$ is a 
product $N'\t\R^k$ in $\C^{m-k}\t\C^k$, where $N'$ is special
Lagrangian in $\C^{m-k}$. Degenerate quadrics with dimension
less than $m-1$ are also excluded.

\subsection{Two trivial constructions of evolution data}
\label{qu42}

Next we consider evolution data not arising from the quadric
construction above. The following two examples are rather trivial
constructions of evolution data, which do not yield interesting
SL $m$-folds in~$\C^m$.

\begin{ex} Let $n\ge 2$, choose any nonzero linear or affine map 
$\chi:\R^n\ra\R^n$, and let $P$ be any integral curve of $\chi$, 
regarded as a vector field in $\R^n$. Then $(P,\chi)$ is a set of 
linear or affine evolution data with $m=2$. Furthermore, every
set of evolution data with $m=2$ comes from this construction.

Thus, using the method of \S\ref{qu3}, one can construct many 
examples of special Lagrangian 2-folds in $\C^2$. But special 
Lagrangian 2-folds in $\C^2$ are equivalent to holomorphic 
curves with respect to an alternative complex structure, and 
so are anyway very easy to construct.
\label{qu4ex4}
\end{ex}

\begin{ex} Let $(P,\chi)$ be a set of evolution data in $\R^n$,
with $P$ an $(m\!-\!1)$-manifold, and let $k\ge 1$. Write
$\R^{n+k}=\R^n\t\R^k$, with coordinates $(x_1,\ldots,x_n$,\allowbreak
$x_{n+1},\ldots,x_{n+k})$. Define
\begin{equation*}
P'=P\t\R^k\quad\text{and}\quad
\chi'=\chi\w{\ts{\pd\over\pd x_{n+1}}}\w\cdots
\w{\ts{\pd\over\pd x_{n+k}}}.
\end{equation*}
Then $(P',\chi')$ is a set of evolution data in $\R^{n+k}$,
with $P'$ an $(m\!+\!k\!-\!1)$-manifold. All SL $(m\!+\!k)$-folds
$N'$ in $\C^{m+k}$ constructed using $(P',\chi')$ split as products
$N\t\R^k$ in $\C^m\t\C^k$, where $N$ is an SL $m$-fold in $\C^m$
constructed using~$(P,\chi)$.
\label{qu4ex5}
\end{ex}

Combining these two examples we can make (uninteresting) examples
of evolution data for any $m,n$ with $2\le m\le n$. In particular,
when $n=m$ we have:

\begin{ex} Let $a,\ldots,f\in\R$ be not all zero, and let $\ga$ be
an integral curve of the vector field $(ax_1+bx_2+e){\pd\over\pd x_1}
+(cx_1+dx_2+f){\pd\over\pd x_2}$ in $\R^2$. Let $m\ge 2$, write
$\R^m=\R^2\t\R^{m-2}$, and define $P=\ga\t\R^{m-2}$ and
\begin{equation*}
\chi=(ax_1\!+\!bx_2\!+\!e){\ts{\pd\over\pd x_1}}\w{\ts{\pd\over\pd x_3}}
\w\cdots\w{\ts{\pd\over\pd x_m}}+(cx_1\!+\!dx_2\!+\!f)
{\ts{\pd\over\pd x_2}}\w\cdots\w{\ts{\pd\over\pd x_m}}.
\end{equation*}
Then $(P,\chi)$ is a set of evolution data with $n=m$. If $e=f=0$
then it is linear, and otherwise affine.
\label{qu4ex6}
\end{ex}

These and the examples of \S\ref{qu41} exhaust the examples with~$m=n$.

\begin{prop} Every set of linear or affine evolution data with\/
$m=n$ is isomorphic either to one of the quadric examples of\/
\S\ref{qu41}, or to one constructed in Example~\ref{qu4ex6}.
\label{qu4prop}
\end{prop}

\begin{proof} Let $(P,\chi)$ be a set of linear or affine evolution
data in $V=\R^m$, with $m=n$. Let $\al$ be a nonzero element of
$\La^mV^*$, and define $\be:V\ra V^*$ by $\be=\al\cdot\chi$, where
`$\,\cdot\,$' is the natural product $\La^mV^*\t\La^{m-1}V\ra V^*$.
Then $\be$ is a linear or affine 1-form on $V$.

The zeros of $\be$ form a distribution $\cal D$ of hyperplanes in $V$
wherever $\be$ is nonzero. The {\it curvature} of $\cal D$ is
$(\d\be)\vert_{\cal D}$. Now clearly $\be\vert_P\equiv 0$, since
$\chi$ is nonzero and tangent to $P$ at each point of $P$, so $\be$
is nonzero along $P$ and ${\cal D}\vert_P=TP$. Therefore $P$ is an
{\it integral submanifold} of $\cal D$, so the curvature of $\cal D$
vanishes along~$P$.

This shows that $\d\be\vert_P\equiv 0$. Clearly, this is equivalent
to $\be\w\d\be$ being zero along $P$, that is, zero in $\La^3V^*$
rather than restricted to $P$. But $\be$ is linear or affine and
$\d\be$ is constant, so $\be\w\d\be$ is linear or affine. As $P$
is not contained in any proper linear or affine subspace of $V$
(as appropriate), we see that $\be\w\d\be$ is zero on all of~$V$.

There are now two possibilities:
\begin{itemize}
\item[(a)] $\d\be=0$, or
\item[(b)] $\d\be=\ga\w\de$ for linearly independent $\ga,\de\in V^*$,
and $\be\in\an{\ga,\de}_{\sst\mathbb R}$ at each point in~$V$.
\end{itemize}
This is because if $\d\be$ is nonzero and not of the form $\ga\w\de$,
then $\be\w\d\be=0$ if and only if $\be=0$, but we know $\be$ is
nonzero on~$P$.

In case (a), we can write $\be=\d Q$ for $Q:V\ra\R$ a quadratic
polynomial, which is homogeneous if $\be$ is linear. Then $Q$ is
constant along $P$ (assuming $P$ connected), so $P$ is a subset
of $P'=\bigl\{v\in V:Q(v)=c\bigr\}$. Thus, case (a) is one of
the quadric examples of \S\ref{qu41}. In case (b), we choose
coordinates $(x_1,\ldots,x_m)$ on $V$ with $\ga=\d x_1$ and
$\de=\d x_2$, and it is then easy to show that we are in the
situation of Example~\ref{qu4ex6}.
\end{proof}

\subsection{Symmetry groups of evolution data}
\label{qu43}

We shall now show that every set of evolution data $(P,\chi)$
has a symmetry group $G$ which is locally transitive on $P$.
For simplicity we work in the linear case; the corresponding
result for affine evolution data may easily be obtained by
replacing linear by affine actions.

\begin{thm} Let\/ $(P,\chi)$ be a set of linear evolution data,
with\/ $P$ a connected, nonsingular $(m\!-\!1)$-submanifold
in $\R^n$. Then there exists a connected Lie subgroup $G$
in $\GL(n,\R)$ with Lie algebra $\g$, such that\/ $P$ is an
open set in a $G$-orbit in $\R^n$, and\/ $\chi$ is $G$-invariant.
Furthermore, there is a natural, surjective, $G$-equivariant
linear map~$L:\La^{m-2}(\R^n)^*\ra\g$.
\label{qu4thm2}
\end{thm}

\begin{proof} Define a linear map $L:\La^{m-1}(\R^n)^*\ra\gl(n,\R)$
by $L(\al)=\chi\cdot\al$, where we regard $\chi$ as an element
of $(\R^n)^*\ot\La^{m-1}\R^n$, and `$\,\cdot\,$' is the natural
contraction $\La^{m-1}\R^n\t\La^{m-2}(\R^n)^*\ra\R^n$, so that
$\chi\cdot\al\in(\R^n)^*\ot\R^n=\gl(n,\R)$. Let $\g$ be the Lie
subalgebra of $\gl(n,\R)$ generated by $\Im L$, so that $L$ maps
$\La^{m-2}(\R^n)^*\ra\g$. Let $G$ be the unique connected Lie
subgroup of $\GL(n,\R)$ with Lie algebra~$\g$.

Regard elements of $\gl(n,\R)$ as linear vector fields on 
$\R^n$. Then at each $p\in P\subset\R^n$ we have $L(\al)\vert_p=
\chi\vert_p\cdot\al$. Since $\chi\vert_p\in\La^{m-1}T_pP$ by
definition, we see that $L(\al)\vert_p\in T_pP$. So the vector
fields $L(\al)$ are tangent to $P$. But the Lie bracket of two
vector fields tangent to $P$ is also tangent to $P$. Hence, as
$\g$ is generated from $\Im L$ by the Lie bracket, every vector
field in $\g$ is tangent to~$P$.

Since $P$ is nonsingular, we have $\chi\vert_p\ne 0$ for all $p\in P$,
by definition. Thus the map $\La^{m-2}(\R^n)^*\ra T_pP$ given by
$\al\mapsto L(\al)\vert_p$ is {\it surjective}. So the vector fields
in $\g$ span $T_pP$ for all $p\in P$. Therefore the action of the Lie
algebra $\g$ on $P$ is {\it locally transitive}. It follows that $P$
is locally isomorphic to an orbit of $G$ in $\R^n$, and as $P$ is
connected, it must be an open set in a $G$-orbit.

Next we prove that $\chi$ is $G$-invariant, which is not quite as
obvious as it looks. Let $1\le i_1<\cdots<i_{m-2}\le n$, set
$\al=\d x_{i_1}\w\cdots\w\d x_{i_{m-2}}$, and define $v=L(\al)$.
We shall show that ${\cal L}_v\chi=0$, where ${\cal L}_v$ is the
Lie derivative. First observe that $v$ is a linear combination of
terms $x_i{\pd\over\pd x_j}$ with $j\ne i_k$ for $k=1,\ldots,m-2$.
It follows easily that ${\cal L}_v\al=0$. But then
\e
0={\cal L}_vv={\cal L}_v(\chi\cdot\al)=
({\cal L}_v\chi)\cdot\al+\chi\cdot({\cal L}_v\al)=
({\cal L}_v\chi)\cdot\al.
\label{qu4eq2}
\e

Now $\chi\vert_P$ is a nonvanishing section of $\La^{m-1}TP$ and $v$
is tangent to $P$, we see that ${\cal L}_v\chi\vert_P=\la\chi\vert_P$
for some smooth function $\la:P\ra\R$. As $({\cal L}_v\chi)\cdot\al=0$
by \eq{qu4eq2}, restricting to $P$ gives $\la\chi\cdot\al=0$ on $P$,
that is, $\la v\equiv 0$ on $P$. Therefore $\la\equiv 0$ or $v\equiv 0$
on $P$. But if $v\equiv 0$ then clearly $\la\equiv 0$. Thus
${\cal L}_v\chi\equiv 0$ on~$P$.

Since $P$ lies in no proper vector subspace of $\R^n$, and ${\cal L}_v\chi$
is linear, this implies that ${\cal L}_v\chi\equiv 0$. This holds whenever
$v=L(\d x_{i_1}\w\cdots\w\d x_{i_{m-2}})$ for $1\le i_1<\cdots<i_{m-2}\le n$.
Such forms are a basis for $\La^{m-2}(\R^n)^*$. So ${\cal L}_v\chi=0$
for all $v\in\Im L$, and therefore for all $v\in\g$. As $G$ is connected,
this shows that $\chi$ is $G$-invariant.

It remains to show that $L:\La^{m-1}(\R^n)^*\ra\g$ is $G$-equivariant
and surjective. The $G$-equivariance is now obvious, as $\chi$ is
$G$-invariant. So $\Im L$ is a $G$-invariant subspace of $\g$, that is,
an {\it ideal} in $\g$. But then $\Im L$ is closed under the Lie bracket.
As $\Im L$ generates $\g$ we have $\g=\Im L$, and $L$ is surjective.
\end{proof}

As an example, consider the linear evolution data $(P,\chi)$ given
in Examples \ref{qu4ex1} and \ref{qu4ex2}. In both cases $G$ is the 
identity component of $\SO(a,m-a)$. In Example \ref{qu4ex1}, each 
connected component of $P$ is an orbit of $G$. In Example \ref{qu4ex2}, 
$P$ is singular at 0, and each component of $P\sm\{0\}$ is an orbit
of~$G$.

Now fix $m=3$. Then $L$ maps $(\R^n)^*\ra\g$. It can be shown that either
\begin{itemize}
\item[(a)] $P$ is contained in no affine hyperplane in $\R^n$,
and $\Ker L=0$; or
\item[(b)] There exists a nonzero linear map $f:\R^n\ra\R$ such that
$P$ is contained in the affine hyperplane $f\equiv 1$ in $\R^n$,
and~$\Ker L=\an{\d f}_{\sst\mathbb R}$.
\end{itemize}
In case (a), $L$ is an isomorphism, so that $\R^n\cong\g^*$. Thus, $P$
is an open set in $G$-orbit in the coadjoint representation $\g^*$ of
$G$, that is, $P$ is locally a {\it coadjoint orbit}. In case (b) we 
will see in \cite[\S 4]{Joyc3} that $(\R^n)^*$ is also a Lie algebra, 
an extension of $\g$ by $\R$, and $P$ is again a coadjoint orbit. Note 
that in case (b) $(P,\chi)$ reduces to a set of affine evolution data 
in~$\R^{n-1}$.

In the sequel to this paper \cite{Joyc3}, we will use these ideas to 
construct a correspondence between sets of evolution data with $m=3$, 
and symplectic 2-manifolds with a transitive, Hamiltonian symmetry group. 
This will enable us to write down several interesting sets of evolution 
data with $m=3$ and $n>3$, and study the corresponding families of SL 
3-folds in~$\C^3$. 

\subsection{Discussion}
\label{qu44}

Let us survey what we know about of sets of evolution data so far.
Evolution data depends on two integers $m,n$ with $2\le m\le n$.
In \S\ref{qu42} we classified all sets of evolution data with
$m=2$ and $m=n$, and constructed some not very interesting examples 
for any $m,n$ with $2\le m\le n$. The ideas of \cite[\S 4]{Joyc3} 
will give us a good picture of the set of all evolution data with 
$m=3$, and could probably be developed into a classification without 
great difficulty.

What we lack at present is an understanding of sets of evolution
data with $3<m<n$. We can state this as:
\medskip

{\bf Problem.} Find and classify examples of sets of evolution
data with $3<m<n$, which do not arise from lower-dimensional 
examples via the product construction of Example~\ref{qu4ex5}.
\medskip

Theorem \ref{qu4thm2} suggests a possible method of constructing
examples. One should start with a likely-looking connected Lie 
group $G$ and a representation $V$ of $G$, and find the 
$(m\!-\!1)$-dimensional orbits $\cal O$ of $G$ in $V$, and then 
look for $G$-invariant elements $\chi$ of $V^*\ot\La^{m-1}V$ which
are tangent to $\cal O$.  Note that if $\chi$ is nonzero and tangent
to $\cal O$ at one point, then it is at every point.

The case $n=m+1$ may also be tractable by a more direct approach.
For instance, affine evolution data with $n=m$, which we understand,
can be interpreted as linear evolution data with $n=m+1$.

Next we discuss the geometric meaning of Theorem \ref{qu4thm2}.
It shows that any set of evolution data $(P,\chi)$ in $\R^n$ has 
a symmetry group $G$, acting on $\R^n$ in a locally transitive
way on $P$ and preserving $\chi$. Take $P$ to be a $G$-orbit
in $\R^n$, so that $G$ acts globally on $P$, rather than just 
locally. Let us ask, how is $G$ related to the special Lagrangian 
$m$-folds $N$ in $\C^m$ constructed from $(P,\chi)$ in~\S\ref{qu3}?

As $N$ is naturally isomorphic to $P\t(-\ep,\ep)$ or $P\t\R$,
and $G$ acts on $P$, there is a natural action of $G$ on $N$.
However, in general this action is {\it not} by automorphisms 
of $\C^m$. That is, $N$ is the image of $\Phi:P\t(-\ep,\ep)\ra\C^m$
and in general there is no $G$-action on $\C^m$ such that $\Phi$ 
is $G$-equivariant.

Nor does $G$ act nontrivially on the set of SL $m$-folds $N$ in
$\C^m$ constructed from $(P,\chi)$. Instead, we should regard
$G$ as acting on the set of {\it parametrizations} $\Phi$ of $N$
constructed in \S\ref{qu3}, so that one SL $m$-fold $N$ will
arise from the construction with many different parametrizations
$\Phi$, related by~$G$.

Here is another way to say this. The maps $\Phi$ were constructed
as solutions of an o.d.e.\ \eq{qu3eq2}, with initial data $\phi_0$ 
in a set ${\cal C}_P$ given in Definitions \ref{qu3def1} and 
\ref{qu3def2}. It turns out that $G$ acts naturally on ${\cal C}_P$, 
and two sets of initial data in the same $G$-orbit in ${\cal C}_P$ 
yield the same SL $m$-fold $N$ in~$\C^m$.

We shall use these ideas to predict the dimension of the family 
${\cal M}_{(P,\chi)}$ of distinct SL $m$-folds $N$ in $\C^m$ 
constructed from $(P,\chi)$ in \S\ref{qu3}. Suppose as above 
that $P$ is a $G$-orbit, and define $G'$ to be the Lie group 
of linear (or affine) automorphisms of $\R^n$ preserving $P$, 
and preserving $\chi$ up to scale. Then $G$ is a subgroup of 
$G'$, but may not be the whole thing. 

For instance, in Example \ref{qu4ex2} $P$ is invariant under 
{\it dilations} ${\bf x}\mapsto t{\bf x}$ in $\R^m$ for $t>0$, 
which do not lie in $G$ for $t\ne 1$, and multiply $\chi$ by 
$t^{m-2}$. In this case $G$ is the identity component of 
$\SO(a,m-a)$, and the identity component of $G'$ is $G\t\R_+$, 
that is, $G$ together with the dilations. In \cite{Joyc3} we 
will give other examples where $G$ needs to be augmented by 
a `dilation' group, which acts in a more complex way on~$\R^n$.

We construct $N$ from the integral curve of an o.d.e.\ in ${\cal C}_P$. 
The set of such curves has dimension $\dim{\cal C}_P-1$. Two curves 
give the same SL $m$-fold $N$ if they are equivalent under the action
of $G'$ on ${\cal C}_P$. Supposing that $G'$ acts locally freely
on ${\cal C}_P$, we guess that~$\dim{\cal M}_{(P,\chi)}=\dim{\cal C}_P
-1-\dim G'$.

In doing this calculation we have factored out the `internal' symmetry
group $G'$ of the construction, which acts on the data used in the
construction, but not on the set of SL $m$-folds we construct. However,
there still remains the `external' symmetry group of automorphisms of
$\C^m$, which is $\SU(m)$ in the linear case (where the origin is a
privileged point) and $\SU(m)\lt\C^m$ in the affine case.

Thus, if generic $m$-folds in ${\cal M}_{(P,\chi)}$ have no continuous 
symmetries, then the moduli space of SL $m$-folds up to automorphisms 
of $\C^m$ has dimension $\dim{\cal M}_{(P,\chi)}\!-\!m^2\!+\!1$ in the 
linear case, and $\dim{\cal M}_{(P,\chi)}\!-\!m^2-\!2m\!+1$ in the affine 
case. This is probably the best measure of the number of `interesting 
parameters' in the construction, once all symmetries are taken into 
account.

\section{Examples from evolving centred quadrics}
\label{qu5}

We will now apply the construction of \S\ref{qu3} to the family 
of sets of linear evolution data $(P,\chi)$ defined using centred 
quadrics in $\R^m$ in Examples \ref{qu4ex1} and \ref{qu4ex2}. In 
\S\ref{qu51} we reduce the problem to an o.d.e.\ in complex functions 
$w_1,\ldots,w_m$ of a real variable $t$, and in \S\ref{qu52} we rewrite 
the o.d.e.\ in terms of functions $u,\th$ and $\th_1,\ldots,\th_m$ of 
$t$. Then in \S\ref{qu53} we solve the equations explicitly, as far as 
we can; the solutions are written in terms of {\it elliptic integrals}.

Section \ref{qu54} considers global properties of the solutions,
and describes the resulting SL $m$-folds in four
different cases. Finally, \S\ref{qu55} considers one particularly
interesting case in which the time evolution may be {\it periodic} 
in $t$, and investigates the conditions for periodicity.

It turns out that in the case that $P$ is a sphere ${\cal S}^{m-1}$ in 
$\R^m$, the SL $m$-folds we construct have already been 
found using a different method by Lawlor \cite{Lawl}, and completed by 
Harvey \cite[p.~139--143]{Harv}. Lawlor used his examples to prove the 
{\it angle conjecture}, a result on when the union of two $m$-planes 
in $\R^n$ is area-minimizing. The other cases of this section can also 
be studied using Lawlor and Harvey's method, and may well be known to 
them, but the author has not found the other cases published anywhere.

Much of this section runs parallel to the construction of 
$\U(1)^{m-2}$-invariant special Lagrangian cones in $\C^m$ in
\cite[\S 7]{Joyc2} and uses the same ideas, because the o.d.e.s 
involved are very similar. However, the geometric interpretations 
are significantly different.

\subsection{Reduction of the problem to an o.d.e.}
\label{qu51}

Let $1\le a\le m$ and $c\in\R$, with $c>0$ if $a=m$, and define 
$P$ and $\chi$ by
\begin{align}
\!\!\!\!\!\!\!\!\!\!\!\!\!\!\!
P&=\bigl\{(x_1,\ldots,x_m)\in\R^m:
x_1^2\!+\!\cdots\!+\!x_a^2\!-\!x_{a+1}^2\!-\!\cdots\!-\!x_m^2\!=\!c\bigr\},
\label{qu5pdef1}\\
\begin{split}
\!\!\!\!\!\!\!\!\!\!\!\!\!\!\!
\chi&=2\sum_{j=1}^a(-1)^{j-1}x_j\,e_1\w\cdots\w e_{j-1}\w 
e_{j+1}\w\cdots\w e_m\\
\!\!\!\!\!\!\!\!\!\!\!\!\!\!\!
&-2\sum_{j=a+1}^m(-1)^{j-1}x_j\,e_1\w\cdots\w e_{j-1}\w 
e_{j+1}\w\cdots\w e_m,
\end{split}
\label{qu5chidef1}
\end{align}
where $e_j={\pd\over\pd x_j}$. Then $(P,\chi)$ is a set of linear 
evolution data. Consider linear maps $\phi:\R^m\ra\C^m$ of the form
\e
\!\!\!\!\!\!\!\!\!\!\!\!\!\!\!\!\!\!\!\!\!\!\!\!\!\!\!\!\!\!\!
\phi:(x_1,\ldots,x_m)\mapsto(w_1x_1,\ldots,w_mx_m)
\quad\text{for $w_1,\ldots,w_m$ in $\C\sm\{0\}$.} 
\!\!\!\!\!\!
\label{qu5eq1}
\e
Then $\phi$ is injective and $\Im\phi$ is a Lagrangian $m$-plane 
in $\C^m$, so that $\phi$ lies in the subset ${\cal C}_P$ of 
$\Hom(\R^m,\C^m)$ given in Definition~\ref{qu3def1}.

We will see that the evolution equation \eq{qu3eq2} for $\phi$
in ${\cal C}_P$ preserves $\phi$ of the form \eq{qu5eq1}. So,
consider a 1-parameter family $\bigl\{\phi_t:t\in(-\ep,\ep)\bigr\}$ 
given by
\e
\phi_t:(x_1,\ldots,x_m)\mapsto\bigl(w_1(t)x_1,\ldots,w_m(t)x_m\bigr),
\label{qu5eq2}
\e
where $w_1,\ldots,w_m$ are differentiable functions from $(-\ep,\ep)$
to $\C\sm\{0\}$. We shall rewrite \eq{qu3eq2} as a first-order o.d.e.\
upon~$w_1,\ldots,w_m$.

Now $(\phi_t)_*(e_j)=w_j{\pd\over\pd z_j}+\bar w_j{\pd\over\pd\bar z_j}$. 
It is convenient to get rid of the $\bar z_j$ term by taking the 
(1,0)-component, giving $(\phi_t)_*(e_j)^{(1,0)}=w_j{\pd\over\pd z_j}$.
In the same way, from \eq{qu5chidef1} the $(m\!-\!1,0)$ component
$(\phi_t)_*(\chi)^{(m-1,0)}$ of $(\phi_t)_*(\chi)$ is
\begin{align*}
&2\sum_{j=1}^a
(-1)^{j\!-\!1}x_jw_1\cdots w_{j\!-\!1}w_{j\!+\!1}\cdots w_m
{\pd\over\pd z_1}\!\w\cdots\w\!{\pd\over\pd z_{j\!-\!1}}\!\w\! 
{\pd\over\pd z_{j\!+\!1}}\!\w\cdots\w\!{\pd\over\pd z_m}\\
-&2\sum_{\!\!\!\!j=a+1\!\!\!\!}^m
(-1)^{j\!-\!1}x_jw_1\cdots w_{j\!-\!1}w_{j\!+\!1}\cdots w_m
{\pd\over\pd z_1}\!\w\cdots\w\!{\pd\over\pd z_{j\!-\!1}}\!\w\! 
{\pd\over\pd z_{j\!+\!1}}\!\w\cdots\w\!{\pd\over\pd z_m}.
\end{align*}

As $\Om$ is an $(m,0)$-tensor, we see that the contraction
of $(\phi_t)_*(\chi)$ with $\Om$ is the same as that of 
$(\phi_t)_*(\chi)^{(m-1,0)}$ with $\Om$. Hence, using the index 
notation for tensors on $\C^m$, we get
\begin{align*}
&(\phi_t)_*(\chi(x))^{a_1\ldots a_{m-1}}\Om_{a_1\ldots a_{m-1}a_m}=\\
&2\sum_{j=1}^ax_jw_1\!\cdots w_{j\!-\!1}w_{j\!+\!1}\!\cdots w_m(\d z_j)_{a_m}
\!-2\sum_{j=a+1}^mx_jw_1\!\cdots w_{j\!-\!1}w_{j\!+\!1}\!\cdots w_m(\d z_j)_{a_m}.
\end{align*}
Hence
\begin{align*}
&(\phi_t)_*(\chi(x))^{a_1\ldots a_{m-1}}
\Om_{a_1\ldots a_{m-1}a_m}g^{a_mb}=\\
&2\sum_{j=1}^ax_jw_1\!\cdots w_{j\!-\!1}w_{j\!+\!1}\!\cdots w_m
\Bigl({\pd\over\pd\bar z_j}\Bigr)^b\!
-2\sum_{j=a+1}^mx_jw_1\!\cdots w_{j\!-\!1}w_{j\!+\!1}\!\cdots w_m
\Bigl({\pd\over\pd\bar z_j}\Bigr)^b.
\end{align*}
Since $(\phi_t)_*(\chi(x))$ and $g$ are real tensors, taking
real parts gives
\begin{align*}
&(\phi_t)_*(\chi(x))^{a_1\ldots a_{m-1}}
(\Re\Om)_{a_1\ldots a_{m-1}a_m}g^{a_mb}=\\
&\sum_{j=1}^ax_j\bar w_1\cdots\bar w_{j\!-\!1}\bar w_{j\!+\!1}\cdots\bar w_m
\Bigl({\pd\over\pd z_j}\Bigr)^b-
\sum_{j=a+1}^mx_j\bar w_1\cdots\bar w_{j\!-\!1}\bar w_{j\!+\!1}\cdots\bar w_m
\Bigl({\pd\over\pd z_j}\Bigr)^b\\
+&\sum_{j=1}^ax_jw_1\cdots w_{j\!-\!1}w_{j\!+\!1}\cdots w_m
\Bigl({\pd\over\pd\bar z_j}\Bigr)^b
-\sum_{j=a+1}^mx_jw_1\cdots w_{j\!-\!1}w_{j\!+\!1}\cdots w_m
\Bigl({\pd\over\pd\bar z_j}\Bigr)^b.
\end{align*}

Now from \eq{qu3eq2} each side of this equation is 
$\bigl({\d\phi_t(x)\over\d t}\bigr)^b$, which satisfies
\begin{equation*}
\Bigl({\d\phi_t(x)\over\d t}\Bigr)^b=
\sum_{j=1}^mx_j{\d w_j\over\d t}\Bigl({\pd\over\pd z_j}\Bigr)^b
+\sum_{j=1}^mx_j{\d\bar w_j\over\d t}\Bigl({\pd\over\pd\bar z_j}\Bigr)^b
\end{equation*}
by \eq{qu5eq2}. Equating coefficients in the last two equations gives
\begin{equation*}
{\d w_j\over\d t}=\begin{cases}
\pha{-}\,\overline{w_1\cdots w_{j-1}w_{j+1}\cdots w_m}, & 
\quad j=1,\ldots,a, \\
-\,\overline{w_1\cdots w_{j-1}w_{j+1}\cdots w_m}, & 
\quad j=a\!+\!1,\ldots,m.
\end{cases}
\end{equation*}
This is the first-order o.d.e.\ upon $w_1,\ldots,w_m$ that we seek.
Applying Theorem \ref{qu3thm2}, we have proved:

\begin{thm} Let\/ $1\le a\le m$ and\/ $c\in\R$, with\/ $c>0$
if\/ $a=m$. Suppose $w_1,\ldots,w_m$ are differentiable functions 
$w_j:(-\ep,\ep)\ra\C\sm\{0\}$ satisfying 
\e
\!\!\!\!\!\!\!\!\!\!\!\!\!\!\!
{\d w_j\over\d t}=\begin{cases}
\pha{-}\,\overline{w_1\cdots w_{j-1} w_{j+1}\cdots w_m},& \quad j=1,\ldots,a, \\
-\,\overline{w_1\cdots w_{j-1} w_{j+1}\cdots w_m},& \quad j=a\!+\!1,\ldots,m.
\end{cases}
\label{qu5eq3}
\e
Define a subset\/ $N$ of\/ $\C^m$ by
\e
\begin{split}
N=\Bigl\{\bigl(w_1(t)x_1&,\ldots,w_m(t)x_m\bigr):t\in(-\ep,\ep),\quad
x_j\in\R,\\
&x_1^2+\cdots+x_a^2-x_{a+1}^2-\cdots-x_m^2=c\Bigr\}.
\end{split}
\label{qu5eq4}
\e
Then $N$ is a special Lagrangian submanifold in~$\C^m$.
\label{qu5thm1}
\end{thm}

Observe that \eq{qu5eq3} agrees with \cite[eq.~(8)]{Joyc2}, with 
$a_j=1$ for $j\le a$ and $a_j=-1$ for $j>a$. Thus, we can follow 
the analysis of \cite[\S 7]{Joyc2} to understand the solutions of 
\eq{qu5eq3}. Furthermore, we showed in \cite[\S 7.6]{Joyc2} that 
\cite[eq.~(8)]{Joyc2} is a {\it completely integrable Hamiltonian 
system}, and the proof also applies to~\eq{qu5eq3}.

\subsection{Rewriting these equations}
\label{qu52}

We now rewrite Theorem \ref{qu5thm1} using different variables. 
If $j\le a$ then 
\eq{qu5eq3} gives
\begin{equation*}
{\d\ms{w_j}\over\d t}=w_j{\d\bar w_j\over\d t}+
\bar w_j{\d w_j\over\d t}=
w_1\cdots w_m+\ov{w_1\cdots w_m}=2\Re(w_1\cdots w_m),
\end{equation*}
and in the same way we get
\e
{\d\ms{w_j}\over\d t}=\begin{cases}
\pha{-}2\Re(w_1\cdots w_m),& \quad j=1,\ldots,a, \\
-2\Re(w_1\cdots w_m),& \quad j=a\!+\!1,\ldots,m.\end{cases}
\label{qu5eq5}
\e

Let $\la\in\R$ be a constant, to be chosen later. Define 
$\al_1,\ldots,\al_m$ by
\e
\al_j=\begin{cases} \ms{w_j(0)}-\la, & \quad j=1,\ldots,a, \\
\ms{w_j(0)}+\la, & \quad j=a+1,\ldots,m, \end{cases}
\label{qu5eq6}
\e
and a function $u:(-\ep,\ep)\ra\R$ by
\begin{equation*}
u(t)=\la+2\int_0^t\Re\bigl(w_1(s)\cdots w_m(s)\bigr)\d s,
\end{equation*}
so that $u(0)=\la$. Then \eq{qu5eq5} gives
\e
\ms{w_j}=\begin{cases} \al_j+u, & \quad j=1,\ldots,a, \\
\al_j-u, & \quad j=a\!+\!1,\ldots,m. \end{cases}
\label{qu5eq7}
\e
Thus we may write
\begin{equation*}
w_j(t)=\begin{cases}{\rm e}^{i\th_j(t)}\sqrt{\al_j+u(t)}, 
& \quad j=1,\ldots,a,\\
{\rm e}^{i\th_j(t)}\sqrt{\al_j-u(t)}, & 
\quad j=a\!+\!1,\ldots,m, \end{cases}
\end{equation*}
for differentiable functions $\th_1,\ldots,\th_m:(-\ep,\ep)\ra\R$. 

Define 
\begin{equation*}
\th=\th_1+\cdots+\th_m \quad\text{and}\quad
Q(u)=\prod_{j=1}^a(\al_j+u)\prod_{j=a+1}^m(\al_j-u).
\end{equation*}
Then we see that 
\begin{equation*}
{\d u\over\d t}=2\Re(w_1\cdots w_m)=2Q(u)^{1/2}\cos\th.
\end{equation*}
Furthermore, expanding out \eq{qu5eq3} shows that
\begin{equation*}
{\d\th_j\over\d t}=\begin{cases}
{\displaystyle -\,{Q(u)^{1/2}\sin\th\over \al_j+u}}, 
& \quad j=1,\ldots,a, \\
{\displaystyle\pha{-}\,{Q(u)^{1/2}\sin\th\over \al_j-u}},
& \quad j=a\!+\!1,\ldots,m.
\end{cases}
\end{equation*}
Summing this equation from $j=1$ to $m$ gives
\begin{equation*}
{\d\th\over\d t}=-\,Q(u)^{1/2}\sin\th
\left(\sum_{j=1}^a{1\over\al_j+u}-\sum_{j=a+1}^m{1\over\al_j-u}\right).
\end{equation*}
Thus, we may rewrite Theorem \ref{qu5thm1} in the following way.

\begin{thm} Let\/ $u$ and\/ $\th_1,\ldots,\th_m$ be differentiable 
functions $(-\ep,\ep)\ra\R$ satisfying
\begin{align}
{\d u\over\d t}&=2Q(u)^{1/2}\cos\th
\label{qu5eq8}\\
\text{and}\quad {\d\th_j\over\d t}&=\begin{cases}
{\displaystyle -\,{Q(u)^{1/2}\sin\th\over \al_j+u}},
& \quad j=1,\ldots,a, \\
{\displaystyle\pha{-}\,{Q(u)^{1/2}\sin\th\over \al_j-u}},
& \quad j=a\!+\!1,\ldots,m,
\end{cases}
\label{qu5eq9}
\end{align}
where $\th=\th_1+\cdots+\th_m$, so that
\e
{\d\th\over\d t}=-\,Q(u)^{1/2}\sin\th
\left(\sum_{j=1}^a{1\over\al_j+u}-\sum_{j=a+1}^m{1\over\al_j-u}\right).
\label{qu5eq10}
\e
Suppose that\/ $\al_j+u>0$ for $j=1,\ldots,a$ and\/ $\al_j-u>0$ for
$j=a\!+\!1,\ldots,m$ and\/ $t\in(-\ep,\ep)$. Define a subset\/ $N$
of\/ $\C^m$ to be
\e
\begin{split}
\!\!\!\!\!\!\!\!\!\!\!\!\!\!\!\!\!\!\!\!\!\!\!\!\!\!\!\!\!\!\!
N=\Bigl\{\bigl(
&x_1{\rm e}^{i\th_1(t)}\sqrt{\al_1+u(t)},
\ldots,x_a{\rm e}^{i\th_a(t)}\sqrt{\al_a+u(t)},
\!\!\!\! \\
\!\!\!\!\!\!\!\!\!\!\!\!
&x_{a+1}{\rm e}^{i\th_{a+1}(t)}\sqrt{\al_{a+1}-u(t)},
\ldots,x_m{\rm e}^{i\th_m(t)}\sqrt{\al_m-u(t)}\,\,\bigr):
\!\!\!\! \\
\!\!\!\!\!\!\!\!\!\!\!\!
&t\in(-\ep,\ep),\;\> x_j\in\R,\;\>
x_1^2\!+\!\cdots\!+\!x_a^2\!-\!x_{a+1}^2\!-\!\cdots\!-\!x_m^2\!=\!c\Bigr\}.
\!\!\!\!
\end{split}
\label{qu5eq11}
\e
Then $N$ is a special Lagrangian submanifold in~$\C^m$.
\label{qu5thm2}
\end{thm}

Now \eq{qu5eq8} and \eq{qu5eq10} give ${\d u\over\d t}$ and 
${\d\th\over\d t}$ as functions of $u$ and $\th$. Dividing 
one by the other gives an expression for ${\d u\over\d\th}$, 
eliminating $t$. Suppose for the moment that $\sin(\th(0))\ne 0$. 
Then separating variables gives
\begin{equation*}
\int_{u(0)}^{u(t)}\left(\sum_{j=1}^a{1\over\al_j+u}
-\sum_{j=a+1}^m{1\over\al_j-u}\right)\d u=
-2\int_{\th(0)}^{\th(t)}\cot\th\,\d\th,
\end{equation*}
which integrates explicitly to
\begin{equation*}
\log Q(u)=-2\log\sin\th+C
\end{equation*}
for all $t\in(-\ep,\ep)$, for some $C\in\R$. Exponentiating 
gives~$Q(u)\sin^2\th\equiv{\rm e}^C>0$. 

If on the other hand $\sin\th(0)=0$ then \eq{qu5eq10} shows that 
$\th$ is constant in $(-\ep,\ep)$, so $Q(u)\sin^2\th\equiv 0$. 
In both cases $Q(u)\sin^2\th$ is constant, so its square root
$Q(u)^{1/2}\sin\th$ is also constant, as it is continuous. Thus 
we have $Q(u)^{1/2}\sin\th\equiv A$ for some~$A\in\R$.

This simplifies \eq{qu5eq9} and \eq{qu5eq10}, as we can replace
the factor $Q(u)^{1/2}\sin\th$ by $A$. Also, from \eq{qu5eq8} we
find that
\begin{equation*}
\Bigl({\d u\over\d t}\Bigr)^2=4Q(u)\cos^2\th=
4\bigl(Q(u)-Q(u)\sin^2\th\bigr)=4\bigl(Q(u)-A^2\bigr).
\end{equation*}
Thus we have proved the following analogue of~\cite[Prop.~7.3]{Joyc2}:

\begin{prop} In the situation of Theorem \ref{qu5thm2} we have
\e
Q(u)^{1/2}\sin\th\equiv A
\label{qu5eq12}
\e
for some $A\in\R$ and all\/ $t\in(-\ep,\ep)$, and\/ 
\eq{qu5eq8}--\eq{qu5eq10} are equivalent to
\begin{align}
\Bigl({\d u\over\d t}\Bigr)^2&=4\bigl(Q(u)-A^2\bigr),
\label{qu5eq13}\\
{\d\th_j\over\d t}&=\begin{cases}
{\displaystyle -\,{A\over\al_j+u}},
& \quad j=1,\ldots,a,\\
{\displaystyle\pha{-}\,{A\over\al_j-u}},
& \quad j=a\!+\!1,\ldots,m,
\end{cases}
\label{qu5eq14}\\
\text{and}\quad 
{\d\th\over\d t}&=-\,A
\left(\sum_{j=1}^a{1\over\al_j+u}-\sum_{j=a+1}^m{1\over\al_j-u}\right).
\label{qu5eq15}
\end{align}
\label{qu5prop1}
\end{prop}

\subsection{Explicit solution using elliptic integrals}
\label{qu53}

Next we will write down the SL $m$-fold $N$ of
Theorem \ref{qu5thm2} in a more simple and explicit way. A nice 
way of doing this is to eliminate $t$, and write everything
instead as a function of $u$. Now ${\d u\over\d t}$ has the same 
sign as $\cos\th$ by \eq{qu5eq8}. Thus, if $\cos\th$ changes sign 
in $(-\ep,\ep)$ then we cannot write $t$ as a function of $u$, but 
if $\cos\th$ has constant sign then we can.

Let us assume that $\th(t)\in(-\pi/2,\pi/2)$ for all $t\in(-\ep,\ep)$, 
so that $\cos\th$ is positive. Then \eq{qu5eq13} gives ${\d u\over\d t}
=2\sqrt{Q(u)-A^2}$, and integrating gives
\begin{equation*}
\int_{u(0)}^{u(t)}{\d u\over 2\sqrt{Q(u)-A^2}}=\int_0^t\d t=t.
\end{equation*}
This defines $u$ implicitly as a function of $t$. From \eq{qu5eq13} 
and \eq{qu5eq14} we get
\begin{equation*}
{\d\th_j\over\d u}=\begin{cases}
{\displaystyle -\,{A\over 2(\al_j+u)\sqrt{Q(u)-A^2}}}, 
& \quad j=1,\ldots,a, \\
{\displaystyle\pha{-}\,{A\over 2(\al_j-u)\sqrt{Q(u)-A^2}}}, 
& \quad j=a\!+\!1,\ldots,m.
\end{cases}
\end{equation*}
Integrating these gives expressions for $\th_j$ in terms of $u$,
and we have proved:

\begin{thm} Suppose $\th(t)\in(-\pi/2,\pi/2)$ for all\/ $t\in(-\ep,\ep)$. 
Then the special Lagrangian $m$-fold\/ $N$ of Theorem \ref{qu5thm2} 
is given explicitly by
\begin{align*}
N=\Bigl\{\bigl(
&x_1{\rm e}^{i\th_1(u)}\sqrt{\al_1+u},
\ldots,x_a{\rm e}^{i\th_a(u)}\sqrt{\al_a+u},\\
&\quad x_{a+1}{\rm e}^{i\th_{a+1}(u)}\sqrt{\al_{a+1}-u},
\ldots,x_m{\rm e}^{i\th_m(u)}\sqrt{\al_m-u}\,\,\bigr):\\
&u\in\bigl(u(-\ep),u(\ep)\bigr),\;\> x_j\in\R,\;\>
x_1^2+\cdots+x_a^2-x_{a+1}^2-\cdots-x_m^2=c\Bigr\},
\end{align*}
where the functions $\th_j(u)$ are given by
\begin{equation*}
\th_j(u)=\begin{cases}
{\displaystyle 
\th_j\bigl(u(0)\bigr)-{A\over 2}\int_{u(0)}^u{\d v\over(\al_j+v)
\sqrt{Q(v)-A^2}}}, & \quad j=1,\ldots,a, \\
{\displaystyle
\th_j\bigl(u(0)\bigr)+{A\over 2}\int_{u(0)}^u{\d v\over(\al_j-v)
\sqrt{Q(v)-A^2}}}, & \quad j=a\!+\!1,\ldots,m.
\end{cases}
\end{equation*}
\label{qu5thm3}
\end{thm}

\subsection{A qualitative description of the solutions}
\label{qu54}

We now describe the SL $m$-folds $N$ in $\C^m$ 
emerging from the construction of Theorem \ref{qu5thm2}, dividing
into four cases, depending on the values of $A$ and~$a$.
\medskip

\noindent{\bf Case (a): $A=0$.} 
\smallskip

\noindent When $A=0$, we see from \eq{qu5eq14} that $\th_1,\ldots,\th_m$ 
are constant with $\th_1+\cdots+\th_m=n\pi$ for some $n\in\Z$, and the 
SL $m$-fold $N$ of \eq{qu5eq11} is a subset of the special 
Lagrangian $m$-plane
\begin{equation*}
\bigl\{(x_1{\rm e}^{i\th_1},\ldots,x_m{\rm e}^{i\th_m}):
x_1,\ldots,x_m\in\R\bigr\}.
\end{equation*}
\medskip

Thus the case $A=0$ is not very interesting. If we replace $\th_1$ 
by $\th_1+\pi$ and $t$ by $-t$ then $A$ changes sign, but the manifold 
$N$ of \eq{qu5eq11} is unchanged. So we may assume in the remaining 
cases that $A>0$. Then it turns out that the equations behave very 
differently depending on whether $a=m$ or $a<m$. We consider the
$a=m$ case first.
\medskip

\noindent{\bf Case (b): $a=m$, $c>0$ and $A>0$.} 
\smallskip

\noindent This case has already been studied by Lawlor \cite{Lawl} 
and Harvey \cite[p.~139--143]{Harv}, using somewhat different 
methods. After some changes of notation, one can show that Harvey 
\cite[Th.~7.78, p.~140]{Harv} is equivalent to the case $a=m$, 
$c>0$ of Theorem \ref{qu5thm3}. Lawlor used his examples to prove 
the {\it angle conjecture}, a result on when the union of two $m$-planes 
in $\R^n$ is area-minimizing. When $\al_1=\cdots=\al_m$, the manifolds
are $\SO(m)$-invariant, and are given in~\cite[\S III.3.B]{HaLa}.

When $m\ge 3$, it can be shown that equation \eq{qu5eq3} admits 
solutions on a bounded open interval $(\ga,\de)$ with $\ga<0<\de$, 
such that $u(t)\ra\iy$ as $t\ra\ga_+$ and $t\ra\de_-$, so that the 
solutions cannot be extended continuously outside $(\ga,\de)$. When 
$m=2$, solutions exist on $\R$, with $u(t)\ra\iy$ as $t\ra\pm\iy$,
so we can put `$\ga=-\iy$' and `$\de=\iy$' in this case. 

The SL $m$-fold $N$ defined using the full solution 
interval $(\ga,\de)$ is a closed, embedded special Lagrangian $m$-fold 
diffeomorphic to ${\cal S}^{m-1}\t\R$. It is the total space of a family 
of ellipsoids $P_t$ in $\C^m$, parametrized by $t$. As $t$ approaches 
$\ga$ or $\de$ these ellipsoids go to infinity in $\C^m$, and also 
become more and more spherical.

At infinity, $N$ is asymptotic to order $r^{1-m}$ to the union 
of two special Lagrangian $m$-planes $\R^m$ in $\C^m$ meeting
at 0, and we can think of $N$ as a {\it connected sum} of two 
copies of $\R^m$. These examples are interesting because they
provide local models for the creation of new SL
$m$-folds in Calabi--Yau $m$-folds as connected sums of other
SL $m$-folds, as in \cite[\S 6--\S 7]{Joyc1} 
when~$m=3$.
\medskip

It remains to consider the cases in which $A>0$ and $1\le a\le m-1$.
Recall that the definition \eq{qu5eq6} of $\al_1,\ldots,\al_m$
depended on an arbitrary constant $\la\in\R$. It is easy to show
that there exists a unique $\la\in\R$ such that $\al_j>0$ for all
$j$, and
\e
\sum_{j=1}^a{1\over\al_j}=\sum_{j=a+1}^m{1\over\al_j}.
\label{qu5eq16}
\e
Let us choose this value of~$\la$.

Since $Q(u)=\ms{w_1}\cdots\ms{w_m}$, and $Q(u)\ge A^2>0$ as $A>0$, 
we have $\ms{w_j}>0$ for all $j$. Thus, from \eq{qu5eq7} we see 
that $u(t)$ is confined to the open interval
\e
\bigl(-\min_{1\le j\le a}\al_j,\min_{a+1\le j\le m}\al_j\bigr)
\label{qu5eq17}
\e
for all $t$ for which the solution exists. It follows from
\eq{qu5eq16} that $Q'(0)=0$. As the roots $-\al_1,\ldots,-\al_a,
\al_{a+1},\ldots,\al_m$ of $Q(u)$ are all real and none lie in
\eq{qu5eq17}, zero is the only turning point of $Q$ in the interval
\eq{qu5eq17}. Thus, $Q$ achieves its maximum in \eq{qu5eq17} at 0, 
and~$Q(0)=\al_1\cdots\al_m$.

But $Q(u)\ge A^2$ by \eq{qu5eq12}. Hence, for all $t$ we have
\e
0<A^2\le Q(u)\le \al_1\cdots\al_m.
\label{qu5eq18}
\e
In particular, this shows that $A\le(\al_1\cdots\al_m)^{1/2}$.
We shall divide into two more cases, depending on whether
$A=(\al_1\cdots\al_m)^{1/2}$ or~$A<(\al_1\cdots\al_m)^{1/2}$.
\medskip

\noindent{\bf Case (c): $1\le a\le m-1$ and $A=(\al_1\cdots\al_m)^{1/2}$.}
\smallskip

\noindent In this case, \eq{qu5eq18} gives $\al_1\cdots\al_m\le 
Q(u)\le\al_1\cdots\al_m$, so $Q(u)\equiv\al_1\cdots\al_m$. It easily
follows that $u\equiv 0$, $\cos\th\equiv 0$ and $\sin\th\equiv 1$,
so that $\th\equiv(2n+\ha)\pi$ for some $n\in\Z$. Equation 
\eq{qu5eq14} then gives
\begin{equation*}
\th_j(t)=\begin{cases}\th_j(0)-At/\al_j & \quad j=1,\ldots,a, \\
\th_j(0)+At/\al_j & \quad j=a\!+\!1,\ldots,m.\end{cases}
\end{equation*}
Thus solutions exist for all $t\in\R$. Define 
\begin{equation*}
a_j=\begin{cases}-A/\al_j & j=1,\ldots,a, \\
\pha{-}A/\al_j & j=a\!+\!1,\ldots,m,\end{cases}
\quad\text{and}\quad y_j=\al_j^{1/2}x_j,\quad j=1,\ldots,m.
\end{equation*}
Then we find that $a_1+\cdots+a_m=0$, and $N$ is given by
\begin{align*}
\Bigl\{\bigl({\rm e}^{i(\th_1(0)+a_1t)}y_1,&\ldots,
{\rm e}^{i(\th_m(0)+a_mt)}y_m\bigr):t\in\R,\quad y_j\in\R,\\
&a_1y_1^2+\cdots+a_my_m^2=-Ac\Bigr\}.
\end{align*}

Now apart from the constant phase factors ${\rm e}^{i\th_j(0)}$, 
this is one of the SL $m$-folds constructed in
\cite[Prop.~9.3]{Joyc2} using the `perpendicular symmetry' idea of 
\cite[\S 9]{Joyc2}, with $n=m$ and $G=\U(1)$ or $\R$. When 
$a_1,\ldots,a_m$ are integers, this example is discussed 
in~\cite[Ex.~9.4]{Joyc2}.
\medskip

\noindent{\bf Case (d): $1\le a\le m-1$ and 
$0<A<(\al_1\cdots\al_m)^{1/2}$.}
\smallskip

\noindent This is very similar to case (c) of \cite[\S 7]{Joyc2}, and
following the proof of \cite[Prop.~7.11]{Joyc2} we can show: 

\begin{prop} Suppose $1\le a\le m-1$, and\/ $\al_1,\ldots,\al_m$
satisfy $\al_j>0$ and\/ \eq{qu5eq16}. Let\/ $u(0)$ and\/ 
$\th_1(0),\ldots,\th_m(0)$ be given, such that\/ $\al_j+u(0)>0$ 
for $j=1,\ldots,a$ and\/ $\al_j-u(0)>0$ for $j=a+1,\ldots,m$, and
\begin{equation*}
0<A=Q(u(0))^{1/2}\sin\th(0)<(\al_1\cdots\al_m)^{1/2}, 
\end{equation*}
where $\th(0)=\th_1(0)+\cdots+\th_m(0)$. Then there exist unique 
solutions $u(t)$, $\th_j(t)$ and\/ $\th(t)$ to equations 
\eq{qu5eq8}--\eq{qu5eq10} of Theorem \ref{qu5thm2} for all\/ 
$t\in\R$, with these values at\/ $t=0$. Furthermore $u$ and\/ $\th$ 
are nonconstant and periodic with period\/ $T>0$, and there exist\/ 
$\be_1,\cdots\be_m\in\R$ with\/ $\be_j<0$ when $j=1,\ldots,a$ and\/
$\be_j>0$ when $j=a\!+\!1,\ldots,m$ and\/ $\be_1+\cdots+\be_m=0$,
such that\/ $\th_j(t+T)=\th_j(t)+\be_j$ for $j=1,\ldots,m$ and 
all\/~$t\in\R$.
\label{qu5prop2}
\end{prop}

Here solutions $u,\th$ to equations \eq{qu5eq8} and \eq{qu5eq10}
are periodic with period $T$, just as in \cite[\S 7.5]{Joyc2}. 
Therefore ${\d\th_j\over\d t}$ is periodic with period $T$ by 
\eq{qu5eq9}, which implies that $\th_j(t+T)=\th_j(t)+\be_j$ for some 
$\be_j\in\R$. But ${\d\th_j\over\d t}<0$ for $j\le a$ and 
${\d\th_j\over\d t}>0$ for $j>a$ by \eq{qu5eq14}, so that 
$\be_j<0$ when $j\le a$ and $\be_j>0$ when $j>a$. Also
$\th=\th_1+\cdots+\th_m$, so that $\th(t+T)=\th(t)+\be_1+\cdots+\be_m$.
As $\th$ is periodic with period $T$, we see that~$\be_1+\cdots+\be_m=0$.

What this means is that when $t$ goes through one cycle of length 
$T$, the complex coordinates $z_1,\ldots,z_m$ don't return to 
their starting points, but instead are taken 
to~${\rm e}^{i\be_1}z_1,\ldots,{\rm e}^{i\be_m}z_m$. 

\subsection{Periodic solutions in case (d)}
\label{qu55}

We have seen that in case (d) above, $u$ and $\th$ are periodic
functions with period $T$, but $\th_1,\ldots,\th_m$ are not periodic,
and satisfy $\th_j(t+T)=\th_j(t)+\be_j$ for $\be_1,\ldots,\be_m$ 
real numbers with $\be_j<0$ if $j\le a$ and $\be_j>0$ if 
$j>a$. But $N$ in \eq{qu5eq11} depends only on 
${\rm e}^{i\th_j}$ rather than on $\th_j$, so that 
$\th_1,\ldots,\th_m$ matter only up to multiples of~$2\pi$.

Thus, if $\be_1,\ldots,\be_m$ are integer multiples of $2\pi$,
then the evolution defining $N$ repeats after time $T$. Actually
it's enough for $\be_j$ to be multiples of $\pi$, as we can change 
the sign of $x_j$ in \eq{qu5eq11}. More generally, if $\be_1,\ldots,\be_m$ 
are rational multiples of $\pi$, then the evolution repeats after
time $nT$, where $n>0$ is the lowest common multiple of the 
denominators of the rational factors.

For our later applications, these periodic solutions are more 
interesting than the non-periodic ones, because they give rise to
closed special Lagrangian $m$-folds in $\C^m$ that can be local 
models for singularities of special Lagrangian $m$-folds in 
Calabi--Yau manifolds. But the non-periodic solutions are not closed 
in $\C^m$, and are not suitable as local models in the same way.

Therefore we will study the dependence of the $\be_j$ upon the 
initial data. It is easy to see that $\be_1,\ldots,\be_m$ depend
only on $a,m$, $\al_1,\ldots,\al_m$ and $A$, and not on $u(0)$ or
$\th_1(0),\ldots,\th_m(0)$. Also, from above the $\al_j$ and $A$ 
satisfy
\e
\!\!\!\!\!\!\!\!\!\!\!\!\!\!\!\!\!\!\!\!\!\!\!\!\!\!\!\!\!\!
\al_j>0,\quad
\sum_{j=1}^a{1\over\al_j}=\sum_{j=a+1}^m{1\over\al_j}
\quad\text{and}\quad
0<A<(\al_1\cdots\al_m)^{1/2}.
\label{qu5eq19}
\e
Given any $\al_1,\ldots,\al_m$ and $A$ satisfying these conditions,
there exists a set of initial data $u(0),\th_1(0),\ldots,\th_m(0)$
with these values. For instance, we can fix $u(0)=0$, and then take
any $\th_1(0),\ldots,\th_m(0)$ such that~$\sin\th(0)=
A(\al_1\cdots\al_m)^{-1/2}$.

Consider what happens to the data when we rescale $N$ by a constant
factor $\ka>0$. Calculation shows that we should replace $u,\th_j,\al_j$ 
and $A$ by $u',\th_j',\al_j'$ and $A'$, where
\begin{equation*}
u'(t)=\ka^2u\bigl(\ka^{m-2}t\bigr),\quad
\th_j'(t)=\th_j\bigl(\ka^{m-2}t\bigr),\quad
\al_j'=\ka^2\al_j\quad\text{and}\quad
A'=\ka^mA.
\end{equation*}
These give new solutions to the equations with period $T'=\ka^{2-m}T$, 
and unchanged values of $\be_1,\ldots,\be_m$. The corresponding SL 
$m$-fold $N'$ is~$\ka N=\{\ka{\bf z}:{\bf z}\in N\}$.

We can now do a parameter count. Since $\be_1+\cdots+\be_m=0$, 
there are only $m-1$ independent $\be_j$. These depend on the $m+1$ 
variables $\al_1,\ldots,\al_m$ and $A$, which satisfy one equation 
$\sum_{j=1}^a{1\over\al_j}=\sum_{j=a+1}^m{1\over\al_j}$. Thus, the 
$\be_j$ can be regarded as $m-1$ functions of $m$ variables.
However, rescaling by $\ka$ leaves the $\be_j$ unchanged, but 
removes one degree of freedom from the $\al_j$ and~$A$. 

Thus, in the initial data $\al_1,\ldots,\al_m$ and $A$ there are
only $m-1$ interesting degrees of freedom, and there are $m-1$
independent $\be_j$ depending on them. The obvious conjecture is
that these two sets of $m-1$ parameters correspond, and that the
map from sets of $\al_j$ and $A$ satisfying \eq{qu5eq19} and
sets of $\be_j$ satisfying $\be_1+\cdots+\be_m=0$ is generically
locally surjective, and locally injective modulo rescaling by 
$\ka>0$ as above. In our next few results we shall show that
this is true.

In the following proposition, modelled on \cite[Prop.~7.13]{Joyc2}, 
we regard the $\al_j$ as fixed, and evaluate the limits of the
$\be_j$ as $A\ra 0$ and $A\ra(\al_1\cdots\al_m)^{1/2}$. For 
simplicity we order the $\al_j$ so that $\al_1\le\cdots\le\al_a$ 
and~$\al_{a+1}\le\cdots\le\al_m$.

\begin{prop} Suppose $1\le a<m$ and\/ $\al_1,\ldots,\al_m>0$ satisfy
\e
\begin{gathered}
\!\!\!\!\!\!\!\!\!\!\!\!\!\!\!\!\!\!\!\!\!\!\!\!\!\!\!\!\!
\al_1=\cdots=\al_k<\al_{k+1}\le\cdots\le\al_a,\\
\!\!\!\!\!\!\!\!\!\!\!\!\!\!\!\!\!\!\!\!\!\!\!\!\!\!\!\!\!
\al_{a+1}\!\le\!\cdots\!\le\!\al_{m-l}\!<\!\al_{m-l+1}\!
=\!\cdots\!=\!\al_m \;\>\text{and}\;\>
\sum_{j=1}^a{1\over\al_j}=\sum_{j=a+1}^m{1\over\al_j}.\!\!\!\!
\end{gathered}
\label{qu5eq20}
\e
Regarding $\al_1,\ldots,\al_m$ as fixed and letting $A$
vary in $\bigl(0,(\al_1\cdots\al_m)^{1/2}\bigr)$, we find that
as $A\ra 0$, we have
\begin{align}
&\be_j\ra\begin{cases} -{\displaystyle{\pi\over k}}, & \quad 1\le j\le k, \\ 
\pha{-}0, & \quad k<j\le m-l, \\ 
\pha{-}{\displaystyle{\pi\over l}}, & \quad m-l<j\le m,\end{cases}
\label{qu5eq21}\\
\intertext{and as $A\ra(\al_1\cdots\al_m)^{1/2}$, we have}
&\be_j\ra\begin{cases} 
-2\pi\al_j^{-1}\Bigl(2\sum_{i=1}^m\al_i^{-2}\Bigr)^{-1/2}, 
& \quad 1\le j\le a, \\
\pha{-}2\pi\al_j^{-1}\Bigl(2\sum_{i=1}^m\al_i^{-2}\Bigr)^{-1/2}, 
& \quad a\!+\!1\le j\le m.\end{cases}
\label{qu5eq22}
\end{align}
\label{qu5prop3}
\end{prop}

\begin{proof} Let $\ga,\de$ be the minimum and maximum values of $u$. 
Then $-\al_1<\ga<0<\de<\al_m$ and $Q(\ga)=Q(\de)=A^2$, and using the 
ideas of \S\ref{qu53} we find that
\e
\!\!\!\!\!\!\!\!\!\!\!\!\!\!\!\!\!\!\!\!\!\!\!\!\!\!\!\!
\be_j=\begin{cases}
{\displaystyle -\int_\ga^\de{\d v\over(\al_j+v)\sqrt{A^{-2}Q(v)-1}}},
& \quad j=1,\ldots,a, \\
{\displaystyle \pha{-}\int_\ga^\de{\d v\over(\al_j-v)\sqrt{A^{-2}Q(v)-1}}},
& \quad j=a\!+\!1,\ldots,m.
\end{cases}
\label{qu5eq23}
\e
As $A\ra 0$ we have $\ga\ra -\al_1$ and $\de\ra\al_m$. Also, the 
factors $(A^{-2}Q(v)-1)^{-1/2}$ in \eq{qu5eq23} tend to zero, 
except near $\ga$ and $\de$. Hence, as $A\ra 0$, the integrands in 
\eq{qu5eq23} get large near $\ga\approx-\al_1$ and $\de\approx\al_m$, 
and very close to zero in between. 

So to understand the $\be_j$ as $A\ra 0$, it is enough to study 
the integrals \eq{qu5eq23} near $\ga$ and $\de$. We shall model 
them at $\ga$. Then near $v=-\al_1$ we have
\begin{equation*}
Q(v)\approx C(v+\al_1)^k,\qquad\text{where}\qquad 
C=\prod_{i=k+1}^a(\al_i-\al_1)\prod_{i=a+1}^m(\al_i+\al_1).
\end{equation*}
Since $A^2=Q(\ga)$ this gives $A^2\approx C(\ga+\al_1)^k$, so 
that~$\ga\approx A^{2/k}C^{-1/k}-\al_1$.

Therefore, when $v\approx\ga$ we have $A^{-2}Q(v)-1\approx 
A^{-2}C(v+\al_1)^k-1$, so when $A$ is small and $j=1,\ldots,k$
we have
\begin{align*}
\int_\ga^0{\d v\over(\al_j+v)\sqrt{A^{-2}Q(v)-1}}&\approx
\int_{A^{2/k}C^{-1/k}-\al_1}^0{\d v\over(\al_1+v)
\sqrt{A^{-2}C(v+\al_1)^k-1}}\\
&\approx \int_0^\iy{2\d w\over k(w^2+1)}={\pi\over k},
\end{align*}
changing variables to $w=\sqrt{A^{-2}C(v+\al_1)^k-1}$, where in
the second line some surprising cancellations happen, and we have 
also approximated the upper limit $\sqrt{A^{-2}C\al_1^k-1}$ by~$\iy$. 

When $k+1\le j\le a$ and $A$ is small we have
\begin{equation*}
\int_\ga^0{\d v\over(\al_j+v)\sqrt{A^{-2}Q(v)-1}}\approx
\int_{A^{2/k}C^{-1/k}-\al_1}^0{\d v\over \al_j
\sqrt{A^{-2}C(v+\al_1)^k-1}}
\approx 0,
\end{equation*}
and similarly when $a+1\le j\le m$ and $A$ is small we have
\begin{equation*}
\int_\ga^0{\d v\over(\al_j-v)\sqrt{A^{-2}Q(v)-1}}\approx 0.
\end{equation*}
So on $[\ga,0]$ the integrals \eq{qu5eq23} are close to $-\pi/k$
for $1\le j\le k$ and 0 for $j>k$ for small $A$. In the same 
way, on $[0,\de]$ the integrals \eq{qu5eq23} are close to $\pi/l$ 
for $m-l<j\le m$ and 0 for $j\le m-l$ for small $A$. This 
proves~\eq{qu5eq21}.

Next consider the behaviour of $\be_j$ as $A\ra(\al_1\cdots\al_m)^{1/2}$. 
When $A$ is close to $(\al_1\cdots\al_m)^{1/2}$, $u$ is small and 
$\sin\th$ close to 1, so $\th$ remains close to $\pi/2$. Write 
$\th={\pi\over 2}+\phi$, for $\phi$ small. Then, setting 
$Q(u)\approx\al_1\cdots\al_m$ and
\begin{equation*}
\cos\th\approx -\phi,\;\>
\sin\th\approx 1 \;\>\text{and}\;\>
\sum_{j=1}^a{1\over\al_j+u}-\sum_{j=a+1}^m{1\over\al_j-u}
\approx -u\sum_{j=1}^m\al_j^{-2},
\end{equation*}
taking only the highest order terms, equations \eq{qu5eq8}
and \eq{qu5eq10} become
\begin{equation*}
{\d u\over\d t}\approx -2(\al_1\cdots\al_m)^{1/2}\phi\quad\text{and}\quad
{\d\phi\over\d t}\approx u(\al_1\cdots\al_m)^{1/2}\sum_{j=1}^m\al_j^{-2},
\end{equation*}
so that $u$ and $\th$ undergo approximately simple harmonic oscillations 
with period $T=2\pi\bigl(2\al_1\cdots\al_m\sum_{j=1}^m\al_j^{-2}
\bigr)^{-1/2}$. Then \eq{qu5eq9} shows that
\begin{equation*}
{\d\th_j\over\d t}\approx \begin{cases}
-\al_j^{-1}(\al_1\cdots\al_m)^{1/2}, & \quad 1\le j\le a,\\
\pha{-}\al_j^{-1}(\al_1\cdots\al_m)^{1/2}, & 
\quad a+1\le j\le m.\end{cases}
\end{equation*}
So $\be_j\approx{\d\th_j\over\d t}T$, as ${\d\th_j\over\d t}$ is 
approximately constant. This proves~\eq{qu5eq22}.
\end{proof}

We can use these limits to show that the map $\boldsymbol{\be}$
from $\al_1,\ldots,\al_m,A$ to $\be_1,\ldots,\be_m$ with
$\be_1+\cdots+\be_m=0$ is generically locally surjective.

\begin{prop} Regard\/ $\boldsymbol{\be}=(\be_1,\ldots,\be_m)$ 
as a function of\/ $\bigl(\al_1,\ldots,\al_m,A\bigr)$. Then 
$\boldsymbol{\be}$ is a real analytic map from $U$ to $V$, where
\begin{gather*}
U=\Bigl\{(\al_1,\ldots,\al_m,A):\al_j>0,\;\>
\sum_{j=1}^a{1\over\al_j}=\!\sum_{j=a+1}^m{1\over\al_j},\;\>
0<A<(\al_1\cdots\al_m)^{1/2}\Bigr\}\\
\text{and}\quad
V=\bigl\{(x_1,\ldots,x_m)\in(-\iy,0)^a\t(0,\iy)^{m-a}:
x_1+\cdots+x_m=0\bigr\}.
\end{gather*}
When $m=2$ we have $\boldsymbol{\be}(u)=(-\pi,\pi)$ for all\/ $u\in U$. 
When $m\ge 3$, the image $\boldsymbol{\be}(U)$ is $(m\!-\!1)$-dimensional, 
and for a dense open subset of\/ $u\in U$ the derivative 
$\d\boldsymbol{\be}\vert_u:\R^{m+1}\ra\bigl\{(x_1,\ldots,x_m)
\in\R^m:x_1+\cdots+x_m=0\bigr\}$ is surjective.
\label{qu5prop4}
\end{prop}

\begin{proof} From \S\ref{qu54} we know that $\be_j<0$ when $j\le a$ 
and $\be_j>0$ when $j>a$, and $\be_1+\cdots+\be_m=0$, so that 
$\boldsymbol{\be}$ does map $U$ to $V$. As $A,\ga$ and $\de$ are 
clearly real analytic functions of the $\al_j$ and $A$, we see from 
\eq{qu5eq23} that $\boldsymbol{\be}$ is real analytic.

When $m=2$ we must have $a=1$, and going back to Theorem 
\ref{qu5thm1} we see that the equations on $w_1,w_2$ are
${\d w_1\over\d t}=\bar w_2$ and ${\d w_2\over\d t}=-\bar w_1$,
which are {\it real linear}, and admit simple harmonic solutions 
with period $2\pi$ for any nonzero initial data. Translating this 
into the notation of Theorem \ref{qu5thm2}, we find that $u$ and 
$\th$ have period $\pi$ and $\be_1=-\pi$, $\be_2=\pi$ for any 
initial data in~$U$.

Now the limits of $\be_1,\ldots,\be_m$ in \eq{qu5eq22} satisfy
$\sum_{j=1}^m\be_j^2=2\pi^2$. Thus, from \eq{qu5eq22} we see
that the closure $\ov{\boldsymbol{\be}(U)}$ contains a nonempty
open subset of the $(m\!-\!2)$-dimensional real hypersurface 
$\sum_{j=1}^mx_j^2=2\pi^2$ in $V$. This implies that 
$\boldsymbol{\be}(U)$ is at least $(m\!-\!2)$-dimensional.

Since $\boldsymbol{\be}$ is real analytic and $U$ is connected,
there are only two possibilities:
\begin{itemize}
\item[(a)] $\boldsymbol{\be}(U)$ is $(m-1)$-dimensional, or
\item[(b)] $\boldsymbol{\be}(U)$ lies in the real hypersurface 
$\sum_{j=1}^mx_j^2=2\pi^2$ in~$V$. 
\end{itemize}
However, when $m\ge 3$ we can use \eq{qu5eq21} to eliminate
possibility (b). For $\ov{\boldsymbol{\be}(U)}$ must contain
the limit in \eq{qu5eq21}, which satisfies $\sum_{j=1}^m\be_j^2
=\pi^2\bigl({1\over k}+{1\over l}\bigr)$. This lies in 
$\sum_{j=1}^mx_j^2=2\pi^2$ only if ${1\over k}+{1\over l}=2$, 
that is, if $k=l=1$, since~$k,l\ge 1$.

Now the ranges of $k$ and $l$ are $k=1,\ldots,a$ and $l=1,\ldots,m\!-\!a$.
When $m=2$ we are forced to take $k=l=a=1$ and we cannot eliminate
possibility (b), as it is actually true. But when $m\ge 3$ we are
always free to choose $k>1$ or $l>1$, so $\ov{\boldsymbol{\be}(U)}$ 
contains a point not on the hypersurface. Thus (b) is false, so (a)
is true, and $\boldsymbol{\be}(U)$ is $(m\!-\!1)$-dimensional. This
shows that $\d\boldsymbol{\be}\vert_u$ must be surjective at some
$u\in U$, and as this is an open condition and $\boldsymbol{\be}$
is real analytic, $\d\boldsymbol{\be}\vert_u$ is surjective for
a dense open subset of~$u\in U$.
\end{proof}

We immediately deduce the following rough analogue 
of~\cite[Cor.~7.14]{Joyc2}.

\begin{cor} In the situation above, we have $\be_1,\ldots,\be_m\in\pi\Q$
for a dense subset of\/ $u$ in~$U$.
\label{qu5cor}
\end{cor}

Here is the main result of this section.

\begin{thm} For each\/ $m\ge 3$ and\/ $1\le a<m$, the construction
above produces a countably infinite collection of\/ $1$-parameter 
families of distinct special Lagrangian $m$-folds $N$ in $\C^m$ 
parametrized by $c\in\R$, given by
\e
\begin{split}
\!\!\!\!\!\!\!\!\!\!\!\!\!\!\!\!\!\!\!\!\!\!\!\!\!\!\!\!\!\!\!
N=\Bigl\{\bigl(
&x_1{\rm e}^{i\th_1(t)}\sqrt{\al_1+u(t)},
\ldots,x_a{\rm e}^{i\th_a(t)}\sqrt{\al_a+u(t)},
\!\!\!\! \\
\!\!\!\!\!\!\!\!\!\!\!\!
&x_{a+1}{\rm e}^{i\th_{a+1}(t)}\sqrt{\al_{a+1}-u(t)},
\ldots,x_m{\rm e}^{i\th_m(t)}\sqrt{\al_m-u(t)}\,\,\bigr):
\!\!\!\! \\
\!\!\!\!\!\!\!\!\!\!\!\!
&t\in\R,\;\> x_j\in\R,\;\>
x_1^2+\cdots+x_a^2-x_{a+1}^2-\cdots-x_m^2=c\Bigr\},
\!\!\!\!
\end{split}
\label{qu5eq24}
\e
such that
\begin{itemize}
\item[{\rm(a)}] if\/ $c>0$ then $N$ is a closed, nonsingular, 
immersed submanifold diffeomorphic to ${\cal S}^{a-1}\!\t\R^{m-a}
\!\t{\cal S}^1$, or to a free quotient of this by~$\Z_2$,
\item[{\rm(b)}] if\/ $c<0$ then $N$ is a closed, nonsingular, 
immersed submanifold diffeomorphic to $\R^a\!\t{\cal S}^{m-a-1}
\!\t{\cal S}^1$, or to a free quotient of this by~$\Z_2$, and
\item[{\rm(c)}] if\/ $c=0$ then $N$ is a closed, immersed cone with 
an isolated singular point at\/ $0$, diffeomorphic to the cone on 
${\cal S}^{a-1}\!\t{\cal S}^{m-a-1}\!\t{\cal S}^1$, or to a cone 
on a free quotient of this by~$\Z_2$.
\end{itemize}
\label{qu5thm4}
\end{thm}

\begin{proof} We saw at the beginning of \S\ref{qu55} that if
$\be_1,\ldots,\be_m\in\pi\Q$ then the time evolution of Theorem
\ref{qu5thm2} exists for all $t$, and is periodic with period 
$nT$ for some $n\ge 1$. But by Corollary \ref{qu5cor} we have 
$\be_1,\ldots,\be_m\in\pi\Q$ for a dense subset of $u\in U$.
Thus, the construction above yields a countable collection of
families of SL $m$-folds, locally parametrized
by $\be_1,\ldots,\be_m\in\pi\Q$ with~$\be_1+\cdots+\be_m=0$.

Choose one of these families, and define $N$ by \eq{qu5eq24}.
Then $N$ is special Lagrangian by Theorem \ref{qu5thm2}. Let
$P$ be the quadric $x_1^2+\cdots+x_a^2-x_{a+1}^2-\cdots-x_m^2=c$
in $\R^m$. Then $N$ is the image of a map $\Phi:P\t\R\ra\C^m$
taking $\bigl((x_1,\ldots,x_m),t\bigr)$ to the point in $\C^m$
defined in \eq{qu5eq24}. As the factors ${\rm e}^{i\th_j(t)}
\sqrt{\al_j\pm u(t)}$ are always nonzero, $\Phi$ is an {\it immersion} 
except when $x_1=\cdots=x_m=0$, which happens only when~$c=0$. 

Thus, $N$ is a nonsingular immersed submanifold when $c\ne 0$, 
and when $c=0$ it has just one singular point 0 as an immersed 
submanifold. Since $\Phi$ is periodic in $t$ with period $nT$, we 
can instead regard $\Phi$ as a map $P\t{\cal S}^1\ra\C^m$, where 
${\cal S}^1=\R/nT\,\Z$. It is also not difficult to see that the
image $N$ of $\Phi$ is {\it closed}, provided $\Phi$ is periodic.

Now $P$ is diffeomorphic to ${\cal S}^{a-1}\!\t\R^{m-a}$ 
when $c>0$, to $\R^a\!\t{\cal S}^{m-a-1}$ when $c<0$, and to the 
cone on ${\cal S}^{a-1}\!\t{\cal S}^{m-a-1}$ when $c=0$. Thus $N$ 
is diffeomorphic under $\Phi$ as an immersed submanifold to 
${\cal S}^{a-1}\!\t\R^{m-a}\!\t{\cal S}^1$ when $c>0$, to
$\R^a\t{\cal S}^{m-a-1}\!\t{\cal S}^1$ when $c<0$, and to 
the cone on ${\cal S}^{a-1}\!\t{\cal S}^{m-a-1}\!\t{\cal S}^1$ 
when~$c=0$.

It remains only to discuss the parts about free quotients by $\Z_2$ 
in (a)--(c). We could have left these bits out, as the result is true 
without them. The point is this: suppose $\be_j=\pi a_j/b$, for 
integers $a_1,\ldots,a_m$ and $b$ with $\hcf(a_1,\ldots,a_m,b)=1$ and 
$b>0$. Then $w_j(t+bT)=(-1)^{a_j}w_j(t)$, so that $w_j$ has period 
$bT$ if $a_j$ is even, and $2bT$ if $a_j$ is odd. Thus $\Phi$ satisfies
\begin{equation*}
\Phi\bigl((x_1,\ldots,x_m),t\bigr)=
\Phi\bigl(
\bigl((-1)^{a_1}x_1,\ldots,(-1)^{a_m}x_m\bigr),t+bT\bigr).
\end{equation*}

Since $P$ is invariant under $(x_1,\ldots,x_m)\mapsto
\bigl((-1)^{a_1}x_1,\ldots,(-1)^{a_m}x_m\bigr)$, the family of 
quadrics making up $N$ has period $bT$. But it doesn't simply
repeat after time $bT$, but also changes the signs of those $x_j$ 
with $a_j$ odd. Let us regard $\Phi$ as mapping $P\t{\cal S}^1\ra\C^m$, 
where ${\cal S}^1=\R/2bT\,\Z$. Then $\Phi$ is generically 2:1, and filters 
through a map $(P\t{\cal S}^1)/\Z_2\ra\C^m$, where the generator of 
$\Z_2$ acts freely on $P\t{\cal S}^1$ by
\begin{equation*}
\bigl((x_1,\ldots,x_m),t+2bT\,\Z\bigr)\mapsto
\bigl(\bigl((-1)^{a_1}x_1,\ldots,(-1)^{a_n}x_m\bigr),t+bT+2bT\,\Z\bigr).
\end{equation*}
This completes the proof.
\end{proof}

This theorem is analogous to \cite[Th.~7.15]{Joyc2}. However, 
\cite[Th.~7.15]{Joyc2} constructs SL 
$T^{m-1}$-cones $N$ with rather large symmetry groups 
$\Sym^0(N)=\U(1)^{m-2}$, but the most generic cones and 
more general SL $m$-folds constructed above
have $\Sym^0(N)=\{1\}$, so they have rather small symmetry
groups. This is not true for all the $m$-folds of Theorem
\ref{qu5thm4}, but only when $\al_1,\ldots,\al_a$ and 
$\al_{a+1},\ldots,\al_m$ are distinct.

Part (c) of the theorem is interesting, as it provides a
large family of singular special Lagrangian cones in $\C^m$
which are good local models for the singularities of special 
Lagrangian $m$-folds in Calabi--Yau $m$-folds. Parts (a) and
(b) are examples of {\it Asymptotically Conical} special
Lagrangian $m$-folds in $\C^m$, and also give local models 
for how the singularities of part (c) can appear as limits 
of families of nonsingular SL $m$-folds in 
Calabi--Yau $m$-folds.

Here is a crude `parameter count' of the number of distinct
families of special Lagrangian $m$-folds produced by this
construction. Locally the families are parametrized by
by $\be_1,\ldots,\be_m\in\pi\Q$ with $\be_1+\cdots+\be_m=0$.
There are unique integers $a_1,\ldots,a_m,b$ with $\be_j=\pi 
a_j/b$, such that $\hcf(a_1,\ldots,a_m,b)=1$ and $b>0$. But
$a_1+\cdots+a_m=0$, so we can discard~$a_m$. 

Observe that the constructions of this section and of \cite[\S 7]{Joyc2}
are strikingly similar in some ways, despite their differences.
The o.d.e.s \eq{qu5eq3} and \cite[eq.~(8)]{Joyc2} behind the two 
constructions are essentially the same. And although the periodicity 
conditions considered above and in \cite[\S 7.5]{Joyc2} are very 
different, the end results are similar, as above we saw that $N$ 
depends on $m$ integers $a_1,\ldots,a_{m-1},b$ with highest common 
factor 1, whereas after \cite[Th.~7.15]{Joyc2} we concluded that $N$ 
depended on $m$ integers $\ti a_1,\ldots,\ti a_{m-1},a$ with 
highest common factor~1. 

The author wonders whether there is some deep connection, or duality,
between the constructions of \cite[\S 7]{Joyc2} and this section, which
explains these similarities. This could be an integrable systems
phenomenon, some kind of `B\"acklund transformation' between the
two constructions which respects the periodicity criteria, or 
something to do with mirror symmetry.

\section{The 3-dimensional case}
\label{qu6}

We now specialize to the case $m=3$ in the situation of \S\ref{qu5}.
The special Lagrangian 3-folds we discuss in this section were
also considered from a different point of view by Bryant 
\cite[\S 3.5]{Brya}. Bryant uses Cartan--K\"ahler theory to study 
special Lagrangian 3-folds $L$ in $\C^3$ whose second fundamental 
form $h$ satisfies certain conditions at every point. 

In effect, Bryant shows \cite[Th.~4]{Brya} that $L$ is one of the 
SL 3-folds of parts (b)--(d) of \S\ref{qu54} if and 
only if $h$ has stabilizer $\Z_2$ in a dense open subset of $L$. His
methods are local. They show that the family of such 3-folds is 
finite-dimensional and compute the dimension, but give less 
information on the global nature of~$L$.

Of the four cases (a)--(d) in \S\ref{qu54}, cases (a), (b) and (c)
are already well understood, so we will concentrate on case (d).
Fixing $m=3$, the two possibilities $a=1$ and $a=2$ in this case 
are exchanged by reversing the order of $z_1,z_2,z_3$ and changing 
the sign of $c$, so without loss of generality we shall choose~$a=1$.

We begin by summarizing the results of \S\ref{qu51} and \S\ref{qu52} 
when $m=3$ and $a=1$. From Theorem \ref{qu5thm1} we obtain

\begin{thm} Suppose $w_1,w_2,w_3:(-\ep,\ep)\ra\C\sm\{0\}$ satisfy
\e
{\d w_1\over\d t}=\overline{w_2w_3},\quad
{\d w_2\over\d t}=-\,\overline{w_3w_1}\quad\text{and}\quad
{\d w_3\over\d t}=-\,\overline{w_1w_2}.
\label{qu6eq1}
\e
Let\/ $c\in\R$, and define a subset\/ $N$ of\/ $\C^3$ to be
\e
\!\!\!\!\!\!\!\!\!\!\!\!\!\!\!\!\!\!\!\!\!\!\!\!\!\!\!\!\!\!\!\!
\Bigl\{\bigl(w_1(t)x_1,w_2(t)x_2,w_3(t)x_3\bigr):t\in(-\ep,\ep),\;\>
x_j\in\R,\;\> x_1^2\!-\!x_2^2\!-\!x_3^2\!=\!c\Bigr\}.
\!\!\!\!\!\!
\label{qu6eq2}
\e
Then $N$ is a special Lagrangian submanifold in~$\C^3$.
\label{qu6thm1}
\end{thm}

Combining Theorem \ref{qu5thm2}, Proposition \ref{qu5prop1} and ideas
from \S\ref{qu54}, we get

\begin{prop} In the situation of Theorem \ref{qu6thm1} the 
functions $w_1,w_2,w_3$ may be written 
\begin{equation*}
w_1={\rm e}^{i\th_1}\sqrt{\al_1+u},\quad
w_2={\rm e}^{i\th_2}\sqrt{\al_2-u}\quad\text{and}\quad
w_3={\rm e}^{i\th_3}\sqrt{\al_3-u},
\end{equation*}
so that
\e
\!\!\!\!\!\!\!\!\!\!\!\!\!\!\!\!\!\!\!\!
\ms{w_1}=\al_1+u,\quad \ms{w_2}=\al_2-u\quad\text{and}\quad
\ms{w_3}=\al_3-u,
\label{qu6eq3}
\e
where $\al_j\in\R$ and\/ $u,\th_1,\th_2,\th_3:(-\ep,\ep)\ra\R$ 
are differentiable functions. Define
\begin{equation*}
Q(u)=(\al_1+u)(\al_2-u)(\al_3-u)\quad\text{and}\quad
\th=\th_1+\th_2+\th_3.
\end{equation*}
Then $Q(u)^{1/2}\sin\th\equiv A$ for some $A\in\R$, and\/
$u$ and\/ $\th_j$ satisfy
\e
\begin{alignedat}{2}
\Bigl({\d u\over\d t}\Bigr)^2&=4\bigl(Q(u)-A^2\bigr),\quad&
{\d\th_1\over\d t}&=-\,{A\over\al_1+u},\\
{\d\th_2\over\d t}&={A\over\al_2-u}\quad&\text{and}\qquad 
{\d\th_3\over\d t}&={A\over\al_3-u}.
\end{alignedat}
\label{qu6eq4}
\e
\label{qu6prop1}
If\/ $A\ne 0$ then $w_j$, $u$ and\/ $\th_j$ exist for all\/ $t$ in 
$\R$, not just in~$(-\ep,\ep)$.
\end{prop}

In the last line, the $w_j$ actually exist for all $t$ even when $A=0$.
But in this case at least one of the $w_j$ will become zero at some time
$t$, and then $\th_j$ is undefined at time $t$, and should be regarded as
jumping discontinuously by $\pm\pi$. Note that as $m=3$, by following the
method of \cite[\S 8.2]{Joyc2} we can solve equation \eq{qu6eq4} explicitly
using the Jacobi elliptic functions. But we will not do this here.

We will explain the 3-dimensional analogue of Theorem \ref{qu5thm4}
in a little more detail. From \S\ref{qu55}, when $m=3$ and $a=1$ 
the SL 3-folds of Theorem \ref{qu5thm4} are locally
parametrized by $\be_1,\be_2,\be_3\in\pi\Q$ with $\be_1<0$, 
$\be_2,\be_3>0$ and $\be_1+\be_2+\be_3=0$. We may write such $\be_j$
uniquely as $\be_j=\pi a_j/b$, where $a_1,a_2,a_3,b\in\Z$ with
$b>0$ and $\hcf(a_1,a_2,a_3,b)=1$. Then the family of quadrics 
making up $N$ has period~$bT$.

However, the functions $w_1,w_2,w_3:\R\ra\C$ satisfy 
$w_j(t+bT)=(-1)^{a_j}w_j(t)$ for $t\in\R$. Thus, if $a_j$ is 
odd then $w_j$ actually has period $2bT$ rather than $bT$. The 
family of quadrics making up $N$ still has period $bT$, because 
the quadric $x_1^2-x_2^2-x_3^2=c$ in $\R^3$ is invariant under 
a change of sign of~$x_j$.

Now in describing the topological type of $N$ in parts (a)--(c) 
of Theorem \ref{qu5thm4}, we allowed the possibility of a free 
quotient by $\Z_2$. When $m=3$ and $a=1$, how this $\Z_2$ acts
depends on whether $a_1$ is even or odd. For instance, when $c>0$ 
the quadric $x_1^2-x_2^2-x_3^2=c$ splits into two connected 
components, with $x_1>0$ and $x_1<0$. Replacing $t$ by $t+bT$
maps $x_j$ to $(-1)^{a_j}x_j$. When $a_1$ is even this map
fixes the two components of the quadric, so that $N$ splits 
into two pieces, but when $a_1$ is odd the two components are 
swapped, so that $N$ comes in only one piece.

We shall state two versions of Theorem \ref{qu5thm4} when $m=3$ 
and $a=1$, for the two cases $a_1$ even and $a_1$ odd. Note that 
the sets of triples $(\be_1,\be_2,\be_3)$ with $\be_j\in\pi\Q$
and $a_1$ even, and with $a_1$ odd, are both dense in the set 
of all $(\be_1,\be_2,\be_3)$, so by the argument of Corollary 
\ref{qu5cor} the sets of initial data with $\be_j\in\pi\Q$ and 
$a_1$ even, and with $\be_j\in\pi\Q$ and $a_1$ odd, are both 
dense in the set of all initial data, and for both cases there 
are a countably infinite number of solutions.

Here is the first version, with $a_1$ even.

\begin{thm} The construction above gives a countably infinite 
collection of\/ $1$-parameter families of distinct special 
Lagrangian $3$-folds $N$ in $\C^3$ parametrized by $c\in\R$, 
given by
\begin{align*}
N=\Bigl\{\bigl(&
x_1{\rm e}^{i\th_1(t)}\sqrt{\al_1+u(t)},
x_2{\rm e}^{i\th_2(t)}\sqrt{\al_2-u(t)},
x_3{\rm e}^{i\th_3(t)}\sqrt{\al_3-u(t)}\,\bigr):\\
&t\in\R,\quad x_j\in\R,\quad
x_1^2-x_2^2-x_3^2=c\Bigr\},
\end{align*}
such that
\begin{itemize}
\item[{\rm(a)}] if\/ $c>0$ then $N$ is the union of two distinct
pieces $N_+$ and\/ $N_-=-N_+$, each of which is a closed, nonsingular, 
immersed submanifold diffeomorphic to~${\cal S}^1\t\R^2$. 
\item[{\rm(b)}] if\/ $c<0$ then $N$ is a closed, nonsingular, 
immersed submanifold diffeomorphic to $T^2\t\R$, with\/ $N=-N$, and
\item[{\rm(c)}] if\/ $c=0$ then $N$ is the union of two distinct
pieces $N_+$ and\/ $N_-=-N_+$, each of which is a closed, immersed 
cone on $T^2$, with an isolated singular point at\/~$0$. 
\end{itemize}
\label{qu6thm2}
\end{thm}

Here is the second version, with $a_1$ odd.

\begin{thm} The construction above gives a countably infinite 
collection of\/ $1$-parameter families of distinct special 
Lagrangian $3$-folds $N$ in $\C^3$ parametrized by $c\in\R$, 
given by
\begin{align*}
N=\Bigl\{\bigl(&
x_1{\rm e}^{i\th_1(t)}\sqrt{\al_1+u(t)},
x_2{\rm e}^{i\th_2(t)}\sqrt{\al_2-u(t)},
x_3{\rm e}^{i\th_3(t)}\sqrt{\al_3-u(t)}\,\bigr):\\
&t\in\R,\quad x_j\in\R,\quad
x_1^2-x_2^2-x_3^2=c\Bigr\},
\end{align*}
such that
\begin{itemize}
\item[{\rm(a)}] if\/ $c>0$ then $N$ is a closed, nonsingular, 
immersed submanifold diffeomorphic to~${\cal S}^1\t\R^2$.
\item[{\rm(b)}] if\/ $c<0$ then $N$ is a closed, nonsingular, 
immersed submanifold diffeomorphic to a free quotient of\/
$T^2\t\R$ by $\Z_2$. It can be thought of as the total space 
of a nontrivial real line bundle over the Klein bottle, and 
has only one infinite end, diffeomorphic to~$T^2\t(0,\iy)$.
\item[{\rm(c)}] if\/ $c=0$ then $N$ is a closed, immersed cone 
on $T^2$, with an isolated singular point at\/~$0$.
\end{itemize}
In all three cases we have~$N=-N$.
\label{qu6thm3}
\end{thm}

In part (c) of these two theorems, the author expects the 
$T^2$-cones to be {\it embedded} in nearly all cases.

\subsection{Conformal parametrization of SL cones}
\label{qu61}

Let us now put $c=0$ in Theorem \ref{qu6thm1}, so that the 3-fold
$N$ of \eq{qu6eq2} is a cone. Define $\Si=N\cap{\cal S}^5$, where
${\cal S}^5$ is the unit sphere in $\C^3$. Then $\Si$ is a {\it 
minimal Legendrian surface} in ${\cal S}^5$, as $N$ is a minimal
Lagrangian 3-fold in~$\C^3$.

We shall write down an explicit {\it conformal parametrization}
$\Phi:\R^2\ra\Si$. Now by \cite[p.~32]{FoWo}, a conformal map
from a Riemann surface to a Riemannian manifold is harmonic
if and only if its image in minimal. Thus, as $\Phi$ is conformal 
and its image $\Si$ is minimal, $\Phi$ is harmonic; and so we have 
constructed an {\it explicit harmonic map} $\Phi:\R^2\ra{\cal S}^5$. 
Such maps are of interest to people who study harmonic maps and 
integrable systems. We begin with a preliminary lemma.

\begin{lem} Let\/ $w_j$ be as in Theorem \ref{qu6thm1} and\/ 
$\al_j$ and\/ $u$ be as in Proposition \ref{qu6prop1}. Then
$\Si=N\cap{\cal S}^5$ may be written
\begin{align*}
\Bigl\{\bigl(w_1(t)&x_1,w_2(t)x_2,w_3(t)x_3\bigr):t\in\R,\quad x_j\in\R,\\
&\al_1x_1^2+\al_2x_2^2+\al_3x_3^2=1,\quad x_1^2-x_2^2-x_3^2=0\bigr\}.
\end{align*}
\label{qu6lem1}
\end{lem}

\begin{proof} A point $(w_1x_1,w_2x_2,w_3x_3)$ in $N$ lies in 
$\Si$ if and only if $\ms{w_1}x_1^2+\ms{w_2}x_2^2+\ms{w_3}x_3^2=1$. 
Substituting in \eq{qu6eq3}, this is equivalent to
\begin{equation*}
(\al_1x_1^2+\al_2x_2^2+\al_3x_3^2)+(x_1^2-x_2^2-x_3^2)=1.
\end{equation*}
But by definition $x_1^2-x_2^2-x_3^2=0$, and 
thus~$\al_1x_1^2+\al_2x_2^2+\al_3x_3^2=1$.
\end{proof}

This shows that $\Si$ is naturally isomorphic to $C\t\R$, where
$C$ is given by
\begin{equation*}
C=\bigl\{(x_1,x_2,x_3)\in\R^3:\al_1x_1^2+\al_2x_2^2+\al_3x_3^2=1,\;\>
x_1^2-x_2^2-x_3^2=0\bigr\}.
\end{equation*}
Since we may assume as in \S\ref{qu54} that $\al_j>0$ for $j=1,2,3$,
it follows that $C$ divides into two connected components $C_+$, 
with $x_1>0$, and $C_-$, with $x_1<0$, each of which is diffeomorphic 
to ${\cal S}^1$. This splitting into $C_\pm$ corresponds to the 
splitting of $N$ into $N_\pm$ in part (c) of Theorem \ref{qu6thm3}. 
There is also a corresponding splitting of $\Si$ into~$\Si_\pm$.

Let us parametrize the circle $C_+$ with a parameter $s$, so that
\begin{equation*}
C_+=\bigl\{\bigl(x_1(s),x_2(s),x_3(s)\bigr):s\in\R\bigr\}.
\end{equation*}
This gives a parametrization $\Phi:\R^2\ra\Si_+$ of $\Si_+$, by
\e
\Phi:(s,t)\mapsto\bigl(w_1(t)x_1(s),w_2(t)x_2(s),w_3(t)x_3(s)\bigr).
\label{qu6eq5}
\e
We shall calculate the conditions upon $x_j(s)$ for $\Phi$ to be
conformal, and solve them.

Since the $x_j(s)$ satisfy $\al_1x_1^2+\al_2x_2^2+\al_3x_3^2=1$ and 
$x_1^2-x_2^2-x_3^2=0$, differentiating with respect to $s$ gives
\begin{equation*}
\al_1x_1\dot x_1+\al_2x_2\dot x_2+\al_3x_3\dot x_3=0\quad\text{and}\quad
x_1\dot x_1-x_2\dot x_2-x_3\dot x_3=0,
\end{equation*}
where `$\dot{\pha{x}}$' is ${\d\over\d s}$. Thus the vector 
$(\dot x_1,\dot x_2,\dot x_3)$ is orthogonal to $(\al_1x_1,
\al_2x_2,\al_3x_3)$ and $(x_1,-x_2,-x_3)$, so it is parallel 
to their vector product. This gives
\e
\begin{gathered}
\dot x_1=\ga(\al_2-\al_3)x_2x_3,\qquad
\dot x_2=-\ga(\al_1+\al_3)x_3x_1\\
\text{and}\qquad \dot x_3=\ga(\al_1+\al_2)x_1x_2,
\end{gathered}
\label{qu6eq6}
\e
for some real nonzero function $\ga(s)$. Also, as $x_1,x_2,x_3$ 
satisfy $\al_1x_1^2+\al_2x_2^2+\al_3x_3^2=1$ and $x_1^2-x_2^2-x_3^2=0$, 
we may write
\e
\!\!\!\!\!\!\!\!\!\!\!\!\!\!\!\!\!\!\!\!\!\!\!\!\!\!\!\!
x_1^2={1+(\al_2-\al_3)v\over \al_1+\al_2},\quad
x_2^2={1-(\al_1+\al_3)v\over\al_1+\al_2}\quad\text{and}\quad
x_3^2=v,
\label{qu6eq7}
\e
for some real function~$v(s)$.

Combining equations \eq{qu6eq1}, \eq{qu6eq5} and \eq{qu6eq6} gives
\begin{align*}
{\pd\Phi\over\pd s}&=\ga\bigl((\al_2-\al_3)w_1x_2x_3,
-(\al_1+\al_3)w_2x_3x_1,(\al_1+\al_2)w_3x_1x_2\bigr)\\
\text{and}\quad
{\pd\Phi\over\pd t}&=\bigl(\,\overline{w_2w_3}\,x_1,
-\,\overline{w_3w_1}\,x_2,-\,\overline{w_1w_2}\,x_3\bigr).
\end{align*}
Thus 
\begin{equation*}
g\Bigl({\pd\Phi\over\pd s},{\pd\Phi\over\pd t}\Bigr)=
\ga\bigl((\al_2-\al_3)+(\al_1+\al_3)-(\al_1+\al_2)\bigr)
\Re(w_1w_2w_3)x_1x_2x_3=0,
\end{equation*}
so that ${{\pd\Phi\over\pd s}}$ and ${{\pd\Phi\over\pd t}}$ are orthogonal.

Using equations \eq{qu6eq3} and \eq{qu6eq7} to write 
$\bms{{\pd\Phi\over\pd s}}$ and $\bms{{\pd\Phi\over\pd t}}$ in terms 
of $u$ and $v$, after a lot of cancellation we find that
\begin{align*}
\Bigl\vert{\pd\Phi\over\pd s}\Bigr\vert^2&=
\ga^2\bigl(\al_3+u+(\al_2-\al_3)(\al_1+\al_3)v\bigr)\\
\text{and}\quad
\Bigl\vert{\pd\Phi\over\pd t}\Bigr\vert^2&=
\al_3+u+(\al_2-\al_3)(\al_1+\al_3)v.
\end{align*}
Note that the coefficients of $uv,v^2$ and $uv^2$ in 
$\bms{{\pd\Phi\over\pd s}}$ and the coefficients of $uv,u^2$ 
and $u^2v$ in $\bms{{\pd\Phi\over\pd t}}$ all vanish. From these
equations, we see that if $\ga^2=1$ then $\bms{{\pd\Phi\over\pd s}}
=\bms{{\pd\Phi\over\pd t}}$, so that $\Phi$ is {\it conformal}.

So let us fix $\ga=1$. Then we seek functions $x_1(s),x_2(s),x_3(s)$
satisfying the o.d.e.\ \eq{qu6eq6} with $\ga=1$, and the restrictions
\eq{qu6eq7}. It turns out that we can solve these equations explicitly
in terms of the {\it Jacobi elliptic functions}, to which we now
give a brief introduction. The following material can be found in 
Chandrasekharan~\cite[Ch.~VII]{Chan}. 

For each $k\in[0,1]$, the Jacobi elliptic functions $\sn(t,k)$, $\cn(t,k)$, 
$\dn(t,k)$ with modulus $k$ are the unique solutions to the o.d.e.s
\begin{align*}
\bigl({\textstyle{\d\over\d t}}\sn(t,k)\bigr)^2&=\bigl(1-\sn^2(t,k)\bigr)
\bigl(1-k^2\sn^2(t,k)\bigr),\\
\bigl({\textstyle{\d\over\d t}}\cn(t,k)\bigr)^2&=\bigl(1-\cn^2(t,k)\bigr)
\bigl(1-k^2+k^2\cn^2(t,k)\bigr),\\
\bigl({\textstyle{\d\over\d t}}\dn(t,k)\bigr)^2&=-\bigl(1-\dn^2(t,k)\bigr)
\bigl(1-k^2-\dn^2(t,k)\bigr),
\end{align*}
with initial conditions
\begin{alignat*}{3}
\sn(0,k)&=0,\;\> & \cn(0,k)&=1,\;\> & \dn(0,k)&=1,\\
{\textstyle{\d\over\d t}}\sn(0,k)&=1,\;\>&
{\textstyle{\d\over\d t}}\cn(0,k)&=0,\;\>&
{\textstyle{\d\over\d t}}\dn(0,k)&=0.
\end{alignat*}
They satisfy the identities
\e
\!\!\!\!\!\!\!\!\!\!\!\!\!\!\!\!\!\!\!\!\!\!\!\!
\sn^2(t,k)+\cn^2(t,k)=1 \;\>\text{and}\;\> k^2\sn^2(t,k)+\dn^2(t,k)=1,
\label{qu6eq8}
\e
and the differential equations
\e
\begin{gathered}
\!\!\!\!\!\!\!\!\!\!\!\!\!\!\!\!\!\!
{\textstyle{\d\over\d t}}\sn(t,k)=\cn(t,k)\dn(t,k),\qquad
{\textstyle{\d\over\d t}}\cn(t,k)=-\sn(t,k)\dn(t,k)\\
\!\!\!\!\!\!\!\!\!\!\!\!\!\!\!\!\!\!
\text{and}\qquad {\textstyle{\d\over\d t}}\dn(t,k)=-k^2\sn(t,k)\cn(t,k).
\end{gathered}
\label{qu6eq9}
\e

Returning to equations \eq{qu6eq6} and \eq{qu6eq7}, suppose 
$\al_2\le\al_3$, and define
\begin{gather*}
x_1=(\al_1+\al_2)^{-1/2}\dn(\mu s,\nu),\quad
x_2=(\al_1+\al_2)^{-1/2}\cn(\mu s,\nu)\\
\text{and}\qquad
x_3=(\al_1+\al_3)^{-1/2}\sn(\mu s,\nu),
\end{gather*}
where
\begin{equation*}
\mu=(\al_1+\al_3)^{1/2}\quad\text{and}\quad
\nu^2={\al_3-\al_2\over\al_1+\al_3}.
\end{equation*}
Then from \eq{qu6eq8} and \eq{qu6eq9}, these $x_j$ satisfy \eq{qu6eq6} 
and \eq{qu6eq7} with $v=(\al_1+\al_3)^{-1}\sn^2(\mu s,\nu)$. Drawing 
the above work together, we have proved:

\begin{thm} In the situation above, define $\Phi:\R^2\ra{\cal S}^5$ by 
\e
\begin{split}
\!\!\!\!\!\!\!\!\!\!\!\!\!\!\!\!\!\!\!\!\!\!\!\!
\Phi&:(s,t)\mapsto\bigl((\al_1+\al_2)^{-1/2}\dn(\mu s,\nu)w_1(t),\\
\!\!\!\!\!\!\!\!\!\!\!\!\!\!\!\!\!\!\!\!\!\!\!\!
&(\al_1+\al_2)^{-1/2}\cn(\mu s,\nu)w_2(t),
(\al_1+\al_3)^{-1/2}\sn(\mu s,\nu)w_3(t)\bigr),
\end{split}
\label{qu6eq10}
\e
where $\mu=(\al_1+\al_3)^{1/2}$, $\nu=(\al_3-\al_2)^{1/2}
(\al_1+\al_3)^{-1/2}$ and\/ ${\cal S}^5$ is the unit sphere in 
$\C^3$. Then $\Phi$ is a conformal, harmonic map.
\label{qu6thm4}
\end{thm}

We made the assumption above that $\al_2\le\al_3$. If $\al_2>\al_3$ then 
we can apply the same method, but swapping over $x_2$ and $x_3$, and 
$\al_2$ and $\al_3$, so that
\begin{gather*}
x_1=(\al_1+\al_3)^{-1/2}\dn(\mu s,\nu),\quad
x_2=(\al_1+\al_2)^{-1/2}\sn(\mu s,\nu)\\
\text{and}\qquad
x_3=(\al_1+\al_3)^{-1/2}\cn(\mu s,\nu),
\end{gather*}
where
\begin{equation*}
\mu=(\al_1+\al_2)^{1/2}\quad\text{and}\quad
\nu^2={\al_2-\al_3\over\al_1+\al_2}.
\end{equation*}
Note also that all of our expressions for $x_j(s)$ depend only on 
the linear combinations $\al_1+\al_2$, $\al_1+\al_3$ and $\al_2-\al_3$ 
of $\al_1,\al_2,\al_3$. This is because the $\al_j$ were defined in 
\eq{qu5eq6} up to an arbitrary constant $\la$, and these combinations
are independent of~$\la$.

\subsection{Relation with harmonic tori in $\CP^2$ and ${\cal S}^5$}
\label{qu62}

Theorem \ref{qu6thm4} constructed a family of {\it explicit conformal 
harmonic maps}\/ $\Phi:\R^2\ra{\cal S}^5$. Furthermore, as the cone on
the image of $\Phi$ is Lagrangian, one can show that if $\pi:{\cal S}^5
\ra\CP^2$ is the Hopf projection then $\pi\circ\Phi$ is conformal and 
harmonic, so we also have a family of explicit conformal harmonic 
maps~$\Psi:\R^2\ra\CP^2$.

Now harmonic maps from Riemann surfaces into spheres and projective
spaces are an {\it integrable system}, and have been intensively
studied in the integrable systems literature. For an introduction
to the subject, see Fordy and Wood \cite{FoWo}, in particular the
articles by Bolton and Woodward \cite[p.~59--82]{FoWo}, McIntosh
\cite[p.~205--220]{FoWo} and Burstall and Pedit~\cite[p.~221--272]{FoWo}.

Therefore our examples can be analyzed from the integrable systems
point of view. We postpone this analysis to the sequel \cite{Joyc6}. 
In \cite[\S 5]{Joyc6} we shall realize the SL cones in $\C^3$
constructed in Theorem \ref{qu6thm1} with $c=0$ as special cases of 
a more general construction of special Lagrangian cones in $\C^3$, 
which involves two commuting o.d.e.s. 

Then in \cite[\S 6]{Joyc6} we work through the integrable systems
framework for the corresponding family of harmonic maps $\Psi:\R^2
\ra\CP^2$, showing that they are generically superconformal of 
finite type, and determining their harmonic sequences, Toda 
solutions, algebras of polynomial Killing fields, and spectral 
curves. From the integrable systems point of view, part (c) of
Theorems \ref{qu6thm2} and \ref{qu6thm3} are interesting because 
they construct large families of {\it superconformal harmonic tori} 
in~$\CP^2$.

\section{Examples from evolving non-centred quadrics}
\label{qu7}

We will now apply the construction of \S\ref{qu3} to the family of 
sets of affine evolution data $(P,\chi)$ defined using non-centred 
quadrics in $\R^m$ in Example \ref{qu4ex3} of \S\ref{qu4}. Our
treatment follows \S\ref{qu5} closely, and so we will leave out 
many of the details.

As in Example \ref{qu4ex3}, let $(m-1)/2\le a\le m-1$, and define 
$P$ and $\chi$ by
\begin{align*}
P&=\bigl\{(x_1,\ldots,x_m)\in\R^m:
x_1^2+\cdots+x_a^2-x_{a+1}^2-\cdots-x_{m-1}^2+2x_m=0\bigr\},\\
\chi&=2(-1)^{m-1}e_1\!\w\!\cdots\!\w\!e_{m-1}\!+
2\sum_{j=1}^a(-1)^{j-1}x_je_1\!\w\!\cdots\!\w\!e_{j-1}\!\w\! 
e_{j+1}\!\w\!\cdots\!\w\!e_m\\
&\qquad\qquad-2\sum_{j=a+1}^{m-1}(-1)^{j-1}x_j\,e_1\w\cdots\w 
e_{j-1}\w e_{j+1}\w\cdots\w e_m,
\end{align*}
where $e_j$ is the vector with $x_j=1$ and $x_k=0$ for $j\ne k$.
Then $P$ is nonsingular in $\R^m$, and $(P,\chi)$ is a set of 
{\it affine evolution data}.

Consider affine maps $\phi:\R^m\ra\C^m$ of the form
\e
\phi:(x_1,\ldots,x_m)\mapsto(w_1x_1,\ldots,w_{m-1}x_{m-1},x_m+\be\,)
\label{qu7eq1}
\e
for $w_1,\ldots,w_{m-1}$ in $\C\sm\{0\}$ and $\be\in\C$. Then 
$\phi$ is injective and $\Im\phi$ is an affine Lagrangian 
$m$-plane in $\C^m$, so that $\phi$ lies in the subset 
${\cal C}_P$ of $\Aff(\R^m,\C^m)$ given in Definition~\ref{qu3def2}.

Then as in \S\ref{qu5}, the evolution equation \eq{qu3eq2} for $\phi$
in ${\cal C}_P$ preserves $\phi$ of the form \eq{qu7eq1}. So, consider 
a 1-parameter family $\bigl\{\phi_t:t\in(-\ep,\ep)\bigr\}$ given by
\begin{equation*}
\phi_t:(x_1,\ldots,x_m)\mapsto\bigl(w_1(t)x_1,\ldots,
w_{m-1}(t)x_{m-1},x_m+\be(t)\bigr),
\end{equation*}
where $w_1,\ldots,w_{m-1}:(-\ep,\ep)\ra\C\sm\{0\}$ and 
$\be:(-\ep,\ep)\ra\C$ are differentiable functions. Following
the method of \S\ref{qu51} one can rewrite \eq{qu3eq2} as a 
first-order o.d.e.\ upon $w_1,\ldots,w_{m-1}$ and $\be$. We
end up with the following analogue of Theorem~\ref{qu5thm1}.

\begin{thm} Let $(m-1)/2\le a\le m-1$. Suppose $w_1,\ldots,w_{m-1}:
(-\ep,\ep)\ra\C\sm\{0\}$ and\/ $\be:(-\ep,\ep)\ra\C\sm\{0\}$ are 
differentiable functions satisfying 
\begin{gather}
\!\!\!\!\!\!\!\!\!\!\!\!\!\!\!\!\!\!\!\!\!\!\!
{\d w_j\over\d t}=
\begin{cases}\pha{-}\,\overline{w_1\cdots w_{j-1} w_{j+1}\cdots w_{m-1}}, 
& \quad j=1,\ldots,a, \\
-\,\overline{w_1\cdots w_{j-1} w_{j+1}\cdots w_{m-1}}, & 
\quad j=a\!+\!1,\ldots,m-1,\end{cases}
\label{qu7eq2}\\
\text{and}\qquad
{\d\be\over\d t}=\overline{w_1\cdots w_{m-1}}. 
\label{qu7eq3}
\end{gather}
Define a subset\/ $N$ of\/ $\C^m$ by
\e
\begin{split}
\!\!\!\!\!\!\!\!\!\!\!\!
N=\Bigl\{&\bigl(w_1(t)x_1,\ldots,w_{m-1}(t)x_{m-1},x_m+\be(t)\bigr):
t\in(-\ep,\ep),\\
& x_j\in\R,\quad x_1^2\!+\!\cdots\!+\!x_a^2\!
-\!x_{a+1}^2\!-\!\cdots\!-\!x_{m-1}^2\!+\!2x_m\!=\!0\Bigr\}.
\end{split}
\label{qu7eq4}
\e
Then $N$ is a special Lagrangian submanifold in~$\C^m$.
\label{qu7thm1}
\end{thm}

Now \eq{qu7eq2} shows that the evolution of $w_1,\ldots,w_{m-1}$
is independent of $\be$. Furthermore, equation \eq{qu7eq2} coincides
with equation \eq{qu5eq3} of Theorem \ref{qu5thm1}, with $m$ replaced by 
$m-1$. Thus, we can use the material of \S\ref{qu52}--\S\ref{qu55}
to write $w_1,\ldots,w_{m-1}$ explicitly in terms of elliptic
integrals, and to describe their global behaviour. 

Having found $w_1,\ldots,w_{m-1}$ as functions of $t$, we can then 
use \eq{qu7eq3} to determine the function $\be$. Thus we can solve 
equations \eq{qu7eq2} and \eq{qu7eq3} in a fairly explicit way, and use 
the solution to describe and understand the SL 
$m$-fold $N$ of~\eq{qu7eq4}.

So, following \S\ref{qu52}, let $\la\in\R$, set $\al_j=\ms{w_j(0)}-\la$ for 
$j=1,\ldots,a$ and $\al_j=\ms{w_j(0)}+\la$ for $j=a\!+\!1,\ldots,m\!-\!1$, 
and define $u:(-\ep,\ep)\ra\R$ by $u(t)=\la+2\int_0^t\Re\bigl(w_1(s)\cdots 
w_{m-1}(s)\bigr)\d s$. Then we have
\begin{equation*}
w_j(t)=\begin{cases}{\rm e}^{i\th_j(t)}\sqrt{\al_j+u(t)}, 
& \quad j=1,\ldots,a,\\
{\rm e}^{i\th_j(t)}\sqrt{\al_j-u(t)}, & 
\quad j=a\!+\!1,\ldots,m\!-\!1, \end{cases}
\end{equation*}
for differentiable functions $\th_1,\ldots,\th_{m-1}:(-\ep,\ep)\ra\R$. Define 
\begin{equation*}
\th=\th_1+\cdots+\th_{m-1} \quad\text{and}\quad
Q(u)=\prod_{j=1}^a(\al_j+u)\prod_{j=a+1}^{m-1}(\al_j-u).
\end{equation*}
Then following \eq{qu5eq8}--\eq{qu5eq10} we find that ${\d u\over\d t}=
2Q(u)^{1/2}\cos\th$, and derive expressions for ${\d\th_j\over\d t}$ 
and~${\d\th\over\d t}$.

As in \eq{qu5eq12} we show that $Q(u)^{1/2}\sin\th\equiv A$ for some 
constant $A\in\R$. Now $w_1\ldots w_{m-1}=Q(u)^{1/2}{\rm e}^{i\th}$. 
Thus equation \eq{qu7eq3} gives 
\begin{equation*}
{\d\be\over\d t}=Q(u)^{1/2}(\cos\th-i\sin\th)=\ha{\d u\over\d t}-iA, 
\end{equation*}
as $Q(u)^{1/2}\cos\th=\ha{\d u\over\d t}$ and $Q(u)^{1/2}\sin\th=A$. 
Integrating this gives
\e
\be(t)=C+\ha u(t)-iAt,
\label{qu7eq5}
\e
where $C=\be(0)-\ha u(0)$. As $\be(0)$ is arbitrary we may as well
fix $C=0$. So we obtain the following analogue of Theorem~\ref{qu5thm2}.

\begin{thm} Let\/ $u$ and\/ $\th_1,\ldots,\th_{m-1}$ be differentiable 
functions $(-\ep,\ep)\ra\R$ satisfying
\begin{align*}
{\d u\over\d t}&=2Q(u)^{1/2}\cos\th\\
\text{and}\quad {\d\th_j\over\d t}&=\begin{cases}
{\displaystyle -\,{Q(u)^{1/2}\sin\th\over \al_j+u}},
& \quad j=1,\ldots,a, \\
{\displaystyle\pha{-}\,{Q(u)^{1/2}\sin\th\over \al_j-u}},
& \quad j=a\!+\!1,\ldots,m\!-\!1,
\end{cases}
\end{align*}
where $\th=\th_1+\cdots+\th_{m-1}$, so that
\begin{equation*}
{\d\th\over\d t}=-\,Q(u)^{1/2}\sin\th
\left(\sum_{j=1}^a{1\over\al_j+u}-\sum_{j=a+1}^{m-1}{1\over\al_j-u}\right).
\end{equation*}
Then $u$ and\/ $\th$ satisfy $Q(u)^{1/2}\sin\th\equiv A$ for some $A\in\R$.
Suppose that\/ $\al_j+u>0$ for $j=1,\ldots,a$ and\/ $\al_j-u>0$ for
$j=a\!+\!1,\ldots,m\!-\!1$ and\/ $t\in(-\ep,\ep)$. Define a subset\/ $N$
of\/ $\C^m$ to be
\begin{align*}
\Bigl\{\bigl(x_1&{\rm e}^{i\th_1(t)}\sqrt{\al_1+u(t)},
\ldots,x_a{\rm e}^{i\th_a(t)}\sqrt{\al_a+u(t)},
x_{a+1}{\rm e}^{i\th_{a+1}(t)}\sqrt{\al_{a+1}-u(t)},\\
&\ldots,x_{m-1}{\rm e}^{i\th_{m-1}(t)}\sqrt{\al_{m-1}-u(t)},
x_m+\ha u(t)-iAt\bigr):\\
&t\in(-\ep,\ep),\;\> x_j\in\R,\;\>
x_1^2+\cdots+x_a^2-x_{a+1}^2-\cdots-x_{m-1}^2+2x_m=0\Bigr\}.
\end{align*}
Then $N$ is a special Lagrangian submanifold in~$\C^m$.
\label{qu7thm2}
\end{thm}

As in \S\ref{qu53}, if we assume that $\th(t)\in(-\pi/2,\pi/2)$ 
for $t\in(-\ep,\ep)$ then $u$ is an increasing function of $t$, and
we can choose to regard everything as a function of $u$ rather
than of $t$. This yields the following analogue of
Theorem~\ref{qu5thm3}:

\begin{thm} Suppose $\th(t)\in(-\pi/2,\pi/2)$ for all\/ $t\in(-\ep,\ep)$. 
Then the special Lagrangian $m$-fold\/ $N$ of Theorem \ref{qu7thm1} 
is given explicitly by
\begin{align*}
N=\Bigl\{\Bigl(
&x_1{\rm e}^{i\th_1(u)}\sqrt{\al_1+u},
\ldots,x_a{\rm e}^{i\th_a(u)}\sqrt{\al_a+u},
x_{a+1}{\rm e}^{i\th_{a+1}(u)}\sqrt{\al_{a+1}-u},\\
&\ldots,x_{m-1}{\rm e}^{i\th_{m-1}(u)}\sqrt{\al_{m-1}-u},
x_m+\ha u-iAt(u)\Bigr):\\
&u\in\bigl(u(-\ep),u(\ep)\bigr),\;\> x_j\in\R,\;\>
x_1^2\!+\!\cdots\!+\!x_a^2\!
-\!x_{a+1}^2\!-\!\cdots\!-\!x_{m-1}^2\!+\!2x_m\!=\!0\Bigr\},
\end{align*}
where the functions $\th_j(u)$ and\/ $t(u)$ are given by
\begin{align*}
&\th_j(u)=\begin{cases}
{\displaystyle 
\th_j(0)-{A\over 2}\int_{u(0)}^u{\d v\over(\al_j+v)\sqrt{Q(v)-A^2}}} 
& \quad j=1,\ldots,a, \\
{\displaystyle
\th_j(0)+{A\over 2}\int_{u(0)}^u{\d v\over(\al_j-v)\sqrt{Q(v)-A^2}}} 
& \quad j=a\!+\!1,\ldots,m\!-\!1,
\end{cases}\\
&\text{and}\qquad t(u)=\int_{u(0)}^u{\d v\over 2\sqrt{Q(v)-A^2}}\,.
\end{align*}
\label{qu7thm3}
\end{thm}

This presentation has the advantage of defining $N$ very explicitly,
but the disadvantage that it is only valid for a certain range of $\th$,
and so of $t$. For understanding the global properties of the solutions
$N$, it is better to keep $t$ as the variable, rather than~$u$.

Next we describe the qualitative behaviour of the solutions,
following the analysis of \S\ref{qu54}. We again divide into 
four cases (a)--(d).
\medskip

\noindent{\bf Case (a): $A=0$.} 
\smallskip

\noindent In this case $N$ is an open subset of a special Lagrangian 
plane $\R^m$ in~$\C^m$.
\medskip

\noindent{\bf Case (b): $a=m-1$ and $A>0$.} 
\smallskip

\noindent When $m\ge 4$, we find that \eq{qu7eq2} and \eq{qu7eq3} 
admit solutions on a bounded open interval $(\ga,\de)$ with 
$\ga<0<\de$, such that $u(t)\ra\iy$ as $t\ra\ga_+$ and $t\ra\de_-$, 
so that the solutions cannot be extended continuously outside 
$(\ga,\de)$. For $m=3$ the solutions exist on $\R$, with
$u(t)\ra\iy$ as $t\ra\pm\iy$, so we can put `$\ga=-\iy$' and
`$\de=\iy$' in this case.

The SL $m$-fold $N$ defined using the full
solution interval $(\ga,\de)$ is a closed, embedded special 
Lagrangian $m$-fold diffeomorphic to $\R^m$, the total space 
of a family of paraboloids $P_t$ in $\C^m$, parametrized by 
$t\in(\ga,\de)$. As $t\ra\ga_+$ and $t\ra\de_-$, these 
paraboloids go to infinity in $\C^m$, and also flatten out,
so that they come to resemble hyperplanes~$\R^{m-1}$.

At infinity, $N$ is asymptotic (in a rather weak sense) to the 
union of two special Lagrangian $m$-planes $\R^m$ in $\C^m$ 
intersecting in $\bigl\{(0,\ldots,0,x_m):x_m\in\R\bigr\}$, a 
copy of $\R$. We should think of these two planes as being
joined when $x_m\in(-\iy,0]$, but separated when~$x_m\in(0,\iy)$.

That is, $N$ is a kind of {\it connected sum} of two special
Lagrangian $m$-planes $\R^m$, but a connected sum performed
along an infinite interval $(-\iy,0]$ rather than a single 
point. Note that $N$ can be regarded as a limiting case of
case (b) of \S\ref{qu54}, in which the two special Lagrangian
$m$-planes degenerate from meeting at a point to meeting at
a line, and at the same time the $x_m$ coordinate of their 
point of intersection goes to~$-\iy$.

These solutions are interesting as local models for singularities 
of SL $m$-folds in Calabi--Yau $m$-folds. When $m=3$ we can solve the 
equations very explicitly, and will do so below.
\medskip

For the two remaining cases with $1\le a\le m-2$ and $A>0$, as in
\S\ref{qu54} we choose the constant $\la$ uniquely such that $\al_j>0$ for 
all $j$ and $\sum_{j=1}^a\al_j^{-1}=\sum_{j=a\!+\!1}^{m\!-\!1}\al_j^{-1}$.
Then $0<A^2\le\al_1\cdots\al_{m-1}$.
\medskip

\noindent{\bf Case (c): $1\le a\le m-2$ and 
$A=(\al_1\cdots\al_{m-1})^{1/2}$.}
\smallskip

\noindent As in \S\ref{qu54}, this is one of the SL 
$m$-folds constructed in \cite[Prop.~9.3]{Joyc2} using the 
`perpendicular symmetry' idea of \cite[\S 9]{Joyc2}, this time with 
$n=m\!-\!1$ and $G=\R$. An example of this with $a=1$ and $m=3$ 
is given in~\cite[Ex.~9.6]{Joyc2}.
\medskip

\noindent{\bf Case (d): $1\le a\le m-2$ and
$0<A<(\al_1\cdots\al_{m-1})^{1/2}$.}
\smallskip

\noindent As in \S\ref{qu54}, in this case solutions exist for all 
$t\in\R$, and $u$ and $\th$ are {\it periodic} in $t$, with period 
$T$. For special values of the initial data we may also arrange for 
$w_1,\ldots,w_{m-1}$ to be periodic with period $nT$ for some~$n\ge 1$.

However, by \eq{qu7eq5} we have $\Im\be(t)=\Im\be(0)-At$, and $A>0$. 
Thus $\be$ is {\it never} periodic, and so the time evolution does 
not repeat itself. So there is no point in following the discussion 
of \S\ref{qu55}. The corresponding SL $m$-folds $N$ 
are embedded submanifolds diffeomorphic to $\R^m$. For various reasons, 
they are not credible as local models for singularities of special 
Lagrangian $m$-folds in Calabi--Yau $m$-folds. 
\medskip

Finally, we set $m=3$. In this case equation \eq{qu7eq2} becomes
a {\it real linear o.d.e.}\ in $w_1$ and $w_2$, and so is
far easier to solve. We consider the cases $a=2$ and $a=1$,
corresponding to cases (b) and (d) above, in the next two
examples.

\begin{ex} Put $m=3$ and $a=2$ in Theorem \ref{qu7thm1}.
Then equations \eq{qu7eq2} and \eq{qu7eq3} become
\e
{\d w_1\over\d t}=\bar w_2,\quad
{\d w_2\over\d t}=\bar w_1\quad\text{and}\quad
{\d\be\over\d t}=\overline{w_1w_2}. 
\label{qu7eq6}
\e
The first two equations have solutions
\begin{equation*}
w_1=C{\rm e}^t+D{\rm e}^{-t}
\quad\text{and}\quad
w_2=\bar C{\rm e}^t-\bar D{\rm e}^{-t},
\end{equation*}
where $C=\ha\bigl(w_1(0)+\ov{w_2(0)}\,\bigr)$ and 
$D=\ha\bigl(w_1(0)-\ov{w_2(0)}\,\bigr)$. Therefore
\begin{equation*}
\overline{w_1w_2}=\ms{C}{\rm e}^{2t}-\ms{D}{\rm e}^{-2t}
+2i\Im(C\bar D),
\end{equation*}
and so integrating the third equation of \eq{qu7eq6} gives
\begin{equation*}
\be(t)=
\ha\ms{C}{\rm e}^{2t}+\ha\ms{D}{\rm e}^{-2t}+2i\Im(C\bar D)t+E,
\end{equation*}
where $E=\be(0)-\ha\ms{C}-\ha\ms{D}$. Thus the special
Lagrangian 3-fold $N$ in $\C^3$ defined in \eq{qu7eq4} is 
given parametrically by
\e
\begin{split}
\!\!\!\!\!\!\!\!\!\!\!\!\!\!\!\!\!\!
\Bigl\{&\bigl(
(C{\rm e}^t+D{\rm e}^{-t})x_1,
(\bar C{\rm e}^t-\bar D{\rm e}^{-t})x_2,-\ha(x_1^2+x_2^2)\\
\!\!\!\!\!\!\!\!\!\!\!\!\!\!\!\!\!\!
&+\ha\ms{C}{\rm e}^{2t}+\ha\ms{D}{\rm e}^{-2t}+2i\Im(C\bar D)t+E\bigr):
x_1,x_2,t\in\R\Bigr\}.
\end{split}
\label{qu7eq7}
\e
Here we have used the equation $x_1^2+x_2^2+2x_3=0$ of \eq{qu7eq4} 
to eliminate~$x_3$. 

Equation \eq{qu7eq7} is a very explicit expression for a special 
Lagrangian 3-fold in $\C^3$. Case (a) above, with $A=0$,
corresponds to $\Im(C\bar D)=0$, and in this case $N$ is a subset 
of an affine special Lagrangian 3-plane $\R^3$ in $\C^3$. If
$\Im(C\bar D)\ne 0$ then $N$ is an embedded submanifold diffeomorphic
to $\R^3$, with coordinates~$(x_1,x_2,t)$.
\label{c3ex6}
\end{ex}

\begin{ex} Put $m=3$ and $a=1$ in Theorem \ref{qu7thm1}.
Then equations \eq{qu7eq2} and \eq{qu7eq3} become
\e
{\d w_1\over\d t}=\bar w_2,\quad
{\d w_2\over\d t}=-\bar w_1\quad\text{and}\quad
{\d\be\over\d t}=\overline{w_1w_2}. 
\label{qu7eq8}
\e
The first two equations have solutions
\begin{equation*}
w_1=C{\rm e}^{it}+D{\rm e}^{-it}
\quad\text{and}\quad
w_2=i\bar D{\rm e}^{it}-i\bar C{\rm e}^{-it},
\end{equation*}
where $C=\ha\bigl(w_1(0)-i\ov{w_2(0)}\bigr)$ and 
$D=\ha\bigl(w_1(0)+i\ov{w_2(0)}\bigr)$. Therefore
\begin{equation*}
\overline{w_1w_2}=iC\bar D{\rm e}^{2it}
-i\bar CD{\rm e}^{-2it}+i\bigl(\ms{C}-\ms{D}\bigr),
\end{equation*}
and so integrating the third equation of \eq{qu7eq8} gives
\begin{equation*}
\be(t)=\ha C\bar D{\rm e}^{2it}+\ha\bar CD{\rm e}^{-2it}
+i\bigl(\ms{C}-\ms{D}\bigr)t+E,
\end{equation*}
where $E=\be(0)-\Re(\bar CD)$. Thus the SL 3-fold 
$N$ in $\C^3$ defined in \eq{qu7eq4} is given parametrically by
\e
\begin{split}
\!\!\!\!\!\!\!\!\!\!\!\!\!\!\!\!\!\!\!\!\!\!\!\!\!\!\!
\Bigl\{&\bigl(
(C{\rm e}^{it}+D{\rm e}^{-it})x_1,
(i\bar D{\rm e}^{it}-i\bar C{\rm e}^{-it})x_2,\ha(x_2^2-x_1^2)\\
\!\!\!\!\!\!\!\!\!\!\!\!\!\!\!\!\!\!\!\!\!\!\!\!\!\!\!
&+\ha C\bar D{\rm e}^{2it}+\ha\bar CD{\rm e}^{-2it}
+i\bigl(\ms{C}-\ms{D}\bigr)t+E\bigr):
x_1,x_2,t\in\R\Bigr\}.\!\!\!\!
\end{split}
\label{qu7eq9}
\e
Here we have used the equation $x_1^2-x_2^2+2x_3=0$ of \eq{qu7eq4} 
to eliminate~$x_3$. 

Case (a) above, with $A=0$, corresponds to $\md{C}=\md{D}$,
and in this case $N$ is a subset of an affine special Lagrangian 
3-plane $\R^3$ in $\C^3$. If $\md{C}\ne\md{D}$ then $N$ is an 
embedded submanifold diffeomorphic to $\R^3$, with coordinates
$(x_1,x_2,t)$. The two cases $C=0$ and $D=0$ are constructed by 
\cite[Prop.~9.3]{Joyc2} with $n=2$, $m=3$ and $G=\R$, as in 
\cite[Ex.~9.6]{Joyc2} and case (c) above, with the symmetry group 
$G$ of $N$ acting by~$(x_1,x_2,t)\mapsto(x_1,x_2,t+c)$.
\label{c3ex7}
\end{ex}

\end{document}